\documentclass[a4paper, 10pt, twoside, notitlepage]{amsart}

\usepackage{amsmath,amscd}
\usepackage{amssymb}
\usepackage{amsthm}
\usepackage{comment}
\usepackage{graphicx, xcolor}

\usepackage{mathrsfs}
\usepackage[ocgcolorlinks, linkcolor=blue]{hyperref}

\usepackage{bm}
\usepackage{bbm}
\usepackage{url}

\usepackage[utf8]{inputenc}
\usepackage{mathtools,amssymb}
\usepackage{esint}
\usepackage{tikz}
\usepackage{dsfont}
\usepackage{relsize}
\usepackage{url}
\usepackage{xcolor}
\usepackage{graphicx}
\usepackage{mathrsfs}
\usepackage[shortlabels]{enumitem}
\usepackage{lineno}
\usepackage{amsmath}
\usepackage{enumitem}
\usepackage{amsthm} 
\usepackage{verbatim}
\usepackage{dsfont}
\numberwithin{equation}{section}

\allowdisplaybreaks

\mathtoolsset{showonlyrefs}

\graphicspath{{images/}}

\newtheorem{theorem}{Theorem}[section]
\newtheorem{lemma}[theorem]{Lemma}
\newtheorem{definition}[theorem]{Definition}

\newtheorem{proposition}[theorem]{Proposition}

\newtheorem{remark}[theorem]{Remark}

{\theoremstyle{remark}
\newtheorem{claim}{Claim}}

\author[P. Zimmermann]{Philipp Zimmermann}
\address{Institute of Mathematics, École Polytechnique Fédérale de Lausanne (EPFL), Lausanne, Switzerland}
\email{philipp.zimmermann@epfl.ch}


\newcommand{\C}{{\mathbb C}}
\newcommand{\R}{{\mathbb R}}

\newcommand{\N}{{\mathbb N}}

\newcommand{\eps}{\varepsilon}

\newcommand{\re}{\text{ Re}} 

\newcommand{\tempered}{\mathscr{S}^{\prime}}

\newcommand{\ifourier}{\mathcal{F}^{-1}}

\newcommand{\distr}{\mathscr{D}^{\prime}}

\newcommand{\weak}{\rightharpoonup}
\newcommand{\weakstar}{\rightharpoonup^\ast}

\title[Weakly viscous nonlocal waves]{Space--time Runge approximation for weakly viscous nonlocal waves and identification of coefficients}
\begin{document}
\maketitle
\begin{abstract}
We study very weak solutions and Runge approximation for weakly viscous
nonlocal wave equations with elastic term $(-\Delta)^su$ and dissipative
term $(-\Delta)^r\partial_tu$, where $0\leq r<s$.
For time-dependent potentials satisfying a smallness condition, we
establish a transposition theory for sources in the dual of the natural
space--time energy space. Under explicit nonlocality and nonresonance assumptions,
we prove Runge approximation in
$L^2(0,T;\widetilde H^s(\Omega))$ and in
$L^2(0,T;\widetilde H^s(\Omega))
\cap H^1_0(0,T;\widetilde H^r(\Omega))$.
As an application, partial exterior Dirichlet-to-Neumann measurements
uniquely determine a bounded, nonnegative, time-dependent coefficient
in the nonlinear velocity damping term $q|\partial_tu|^{p-2}\partial_tu$
for $0<r<s$ and $2<p\leq 2+4r/n$.
The proof establishes the small-amplitude response by energy estimates,
weighted $L^p$ monotonicity convergence, and a linear Green identity,
then applies adjoint Runge approximation followed by a forward velocity
approximation in $L^2$. No smallness or positive lower bound is imposed
on the nonlinear damping coefficient.
\end{abstract}

	\tableofcontents

\section{Introduction}
\label{intro:section}
Exterior measurements for nonlocal wave equations contain information about
both the spatial behavior of solutions and their time evolution. Runge
approximation makes this information accessible by approximating prescribed
interior profiles with solutions generated from the exterior. For inverse
problems involving nonlinear functions of the velocity, approximation of the
displacement alone is insufficient: one needs a topology that also controls
its time derivative. The purpose of this article is to establish such a
space--time approximation property for weakly viscous nonlocal waves and to
apply it to the recovery of a nonlinear damping coefficient.

Let $\Omega\subset\R^n$ be a bounded Lipschitz domain,
$\Omega_e=\R^n\setminus\overline\Omega$, and $T>0$. We consider the linear equation\footnote{For a set
$A\subset\R^n$, we write $A_T=A\times(0,T)$.}
\begin{equation}\label{intro:linear-problem}
\begin{cases}
 (\partial_t^2+(-\Delta)^r\partial_t+(-\Delta)^s+q)u=f
 &\text{in }\Omega_T,\\
 u=\varphi&\text{in }(\Omega_e)_T,\\
 u(0)=u_0,\quad\partial_tu(0)=u_1&\text{in }\Omega,
\end{cases}
\end{equation}
where $0\leq r<s<\infty$ and $(-\Delta)^\rho$ is the Fourier multiplier
with symbol $|\xi|^{2\rho}$. Starting from the finite energy theory in
\cite{zimmermann-paperI}, we construct very weak solutions for sources in
the dual of a space--time energy space. This supplies the adjoint theory
needed for stronger Runge approximation. We then study the nonlinear
problem obtained by replacing $qu$ with $q(x,t)|\partial_tu|^{p-2}\partial_tu$.
The coefficient $q$ in this nonlinear problem plays a different role from
the linear potential in \eqref{intro:linear-problem}.

\subsection{Weakly viscous waves and space--time approximation}
The inequality $r<s$ means that the spatial order of the dissipation is
smaller than that of the elastic term. The finite energy theory controls
$u$ in $C([0,T];\widetilde H^s(\Omega))$ and $\partial_tu$ in
$C([0,T];L^2(\Omega))\cap L^2(0,T;\widetilde H^r(\Omega))$ for homogeneous
exterior data. The resulting space--time space is
\[
 X^1_{s,r}(\Omega_T)
 =L^2(0,T;\widetilde H^s(\Omega))
   \cap H^1(0,T;\widetilde H^r(\Omega)).
\]
Here $\widetilde H^\rho(\Omega)$ denotes the closure of
$C_c^\infty(\Omega)$ in $H^\rho(\R^n)$, and $H^{-\rho}(\Omega)$ its
antidual. We use $X^1_{s,r,0}$ and $X^1_{s,r,T}$ for the subspaces with
vanishing displacement trace at $0$ and $T$, respectively. In particular,
the notation $H^1_0$ in this article imposes only the initial trace.

The Runge theorem of \cite{zimmermann-paperI} gives approximation in
$L^2(0,T;\widetilde H^r(\Omega))$, and hence in $L^2(\Omega_T)$.
To approximate in $X^1_{s,r,0}(\Omega_T)$, Hahn--Banach duality requires
sources that can act on time derivatives and traces. These sources are
not generally represented by functions in a spatial Sobolev dual:
\[
 (X^1_{s,r}(\Omega_T))'
 =L^2(0,T;H^{-s}(\Omega))
  +\bigl(H^1(0,T;\widetilde H^r(\Omega))\bigr)'.
\]
We therefore develop a transposition theory for this larger source class.
The construction combines approximation of the rough data with parabolic
regularization of the finite energy problems. An integrated-in-time test
function yields the required estimates under a smallness condition on
the linear potential.

When the elastic order is an integer, the exterior orthogonality
identity initially leaves a time-independent mode in both Runge
arguments. For sources in $L^2(0,T;H^{-s}(\Omega))$, we eliminate
this mode using the transposition identity and the zero initial
conditions. For general energy-dual sources, we use adjoint test
functions concentrated near the initial time. These functions
converge to zero in the space--time energy space while retaining
a prescribed initial velocity. Continuity of the source functional
then forces the stationary mode to vanish. This argument also
applies to $L^2(0,T;H^{-s}(\Omega))$ sources.
Consequently, spatial-energy and full space--time Runge
approximation hold in the same three order regimes described below.

\subsection{Main results}
We first specify the potential class used in the linear theory.
The class $\mathcal Q_s$ consists of real-valued potentials
$q\in L^\infty(0,T;L^\ell(\Omega))$, for some $2\leq\ell\leq\infty$
satisfying
\[
\begin{cases}
 n/s\leq\ell\leq\infty,&2s<n,\\
 2<\ell\leq\infty,&2s=n,\\
 2\leq\ell\leq\infty,&2s>n,
\end{cases}
\]
that are distributionally continuous in time, meaning that the map
\[
    t\mapsto\int_\Omega q(x,t)\zeta(x)\,dx
\]
is continuous on $[0,T]$ for every $\zeta\in C_c^\infty(\Omega)$.
Our very weak theory requires
\begin{equation}\label{intro:smallness}
 \|q\|_{L^1(0,T;L^\ell(\Omega))}\leq(8C_*)^{-1},
\end{equation}
where $C_*>0$ is a numerical constant depending only on $n,r,s,\ell$, and $\Omega$.

\begin{theorem}[Very weak well-posedness]\label{intro:very-weak}
Let $0\leq r<s$ and let $q\in\mathcal Q_s$ satisfy
\eqref{intro:smallness}. Suppose that
\[
 (u_0,u_1)\in
 \begin{cases}
 \widetilde H^{2r-s}(\Omega)\times H^{-s}(\Omega),&s<2r,\\
 L^2(\Omega)\times H^{-s}(\Omega),&2r\leq s,
 \end{cases}
 \qquad f\in(X^1_{s,r}(\Omega_T))'.
\]
Then \eqref{intro:linear-problem} with $\varphi=0$ has a unique very
weak solution
\[
 u\in L^\infty(0,T;\widetilde L^2(\Omega))
       \cap L^2(0,T;\widetilde H^r(\Omega)),
\]
in the transposition sense of
Definition~\ref{F:def: very weak sols for energy sources}.
The solution depends continuously on the data in the indicated spaces.
\end{theorem}

This is Theorem~\ref{F:thm: very weak sol with rough sources}.
We first prove a preliminary result for sources in the dual space
$L^2(0,T;H^{-s}(\Omega))$ and then pass to general energy-dual sources.
The initial conditions in this theorem are encoded in the transposition
identity; pointwise velocity traces are not asserted for very weak solutions.

Let $W\subset\Omega_e$ be nonempty and open. For exterior data
$\varphi\in C_c^\infty(W_T)$, let $u_\varphi$ solve
\eqref{intro:linear-problem} with $f=u_0=u_1=0$, and define
\[
 \mathscr R_W=\{u_\varphi-\varphi:\varphi\in C_c^\infty(W_T)\}.
\]
We distinguish the following regimes, always with $0\leq r<s$:
\begin{enumerate}[(i)]
\item\label{intro:orders-one} $s\in\N$ and $r\notin\N_0$;
\item\label{intro:orders-two} $s\notin\N$ and $r\in\N_0$;
\item\label{intro:orders-three} $n\geq2$, $r,s\notin\N_0$, and
\[
 \begin{cases}
 s-r\notin\N,&n\text{ even},\\
 s-r\notin\tfrac12\N,&n\text{ odd}.
 \end{cases}
\]
\end{enumerate}

\begin{theorem}[Spatial-energy and space--time Runge approximation]
\label{intro:runge}
Let $q\in\mathcal Q_s$ satisfy \eqref{intro:smallness}.
Under any of the regimes \ref{intro:orders-one}--\ref{intro:orders-three},
$\mathscr R_W$ is dense in $L^2(0,T;\widetilde H^s(\Omega))$ and in $X^1_{s,r,0}(\Omega_T)$.
\end{theorem}

These assertions are Theorems~\ref{R:thm: Runge approx in Hs}
and~\ref{R:thm: Runge approx in X}. Their proofs combine the very weak
theory with unique continuation and the entanglement principle for
fractional Laplacians. Time mollification permits these spatial tools
to be applied to the low-regularity transposition solutions. In the
space--time argument, density of suitable adjoint test functions is
also needed to recover the entire annihilating functional, including
its trace contributions. Time reversal gives the corresponding adjoint
approximation in $X^1_{s,r,T}(\Omega_T)$.

Our inverse application concerns
\begin{equation}\label{intro:nonlinear-problem}
\begin{cases}
 (\partial_t^2+(-\Delta)^r\partial_t+(-\Delta)^s)u
   +q(x,t)|\partial_tu|^{p-2}\partial_tu=0&\text{in }\Omega_T,\\
 u=\varphi&\text{in }(\Omega_e)_T,\\
 u(0)=\partial_tu(0)=0&\text{in }\Omega.
\end{cases}
\end{equation}
For this problem, $q\in L^\infty(\Omega_T)$ is nonnegative. The
Dirichlet-to-Neumann map is the distributional exterior flux
\[
 \Lambda_q(\varphi)
 =\left.\bigl((-\Delta)^r\partial_t+(-\Delta)^s\bigr)u\right|_{(\Omega_e)_T}.
\]

\begin{theorem}[Recovery of nonlinear velocity damping]
\label{intro:nonlinear-uniqueness}
Let $0<r<s$ satisfy regime \ref{intro:orders-one}, \ref{intro:orders-two}, or
\ref{intro:orders-three}, and suppose that
\[
 2<p\leq 2+4r/n.
\]
Let $q_1,q_2\in L^\infty(\Omega_T)$ be nonnegative, and let
$W_1,W_2\subset\Omega_e$ be nonempty open sets. If, for every real-valued
$\varphi\in C_c^\infty((W_1)_T)$ and all sufficiently small $\delta>0$,
\[
 \Lambda_{q_1}(\delta\varphi)|_{(W_2)_T}
 =\Lambda_{q_2}(\delta\varphi)|_{(W_2)_T},
\]
then $q_1=q_2$ almost everywhere in $\Omega_T$.
The smallness threshold for $\delta$ may depend on $\varphi$.
\end{theorem}

This is Theorem~\ref{thm:nld-uniqueness}. There is no smallness assumption
on the nonlinear coefficients $q_j$, since the Runge theorem is applied to
the linear equation with zero potential. The sets $W_1$ and $W_2$ need
not be disjoint. Nonlinear forward solvability follows by parabolic
regularization, a local contraction argument, and continuation using
the dissipative energy estimate.

The first step in recovery is the small-amplitude response formula of
Proposition~\ref{prop:nld-response}. If $v_\varphi$ is the linear forward
solution and $w_\zeta$ the adjoint solution with zero terminal data, then
\[
\begin{split}
 &\lim_{\delta\downarrow0}\delta^{-(p-1)}
 \langle\Lambda_q(\delta\varphi)-\delta\Lambda_0(\varphi),\zeta\rangle\\
 &\qquad=\int_{\Omega_T}q|\partial_tv_\varphi|^{p-2}
                  \partial_tv_\varphi\,w_\zeta\,dx\,dt.
\end{split}
\]
We prove this identity using forward energy estimates, strong convergence
of the rescaled velocities in $L^p(\Omega_T,q\,dx\,dt)$, and the linear
source-to-boundary Green identity. No higher-order differentiability of
the nonlinear solution map is needed, and the coefficient may vanish.
Adjoint Runge approximation first yields
\[
    (q_1-q_2)|\partial_t v_\varphi|^{p-2}
    \partial_t v_\varphi=0
    \quad\text{a.e. in }\Omega_T
\]
for each fixed exterior datum $\varphi$.
Forward space--time Runge approximation then provides velocity
convergence in $L^2(\Omega_T)$ and, along a subsequence,
almost everywhere.
Continuity of the nonlinearity identifies the coefficient
throughout the range $2<p\leq2+4r/n$.

\subsection{Relation to inverse problems for nonlocal waves}
The relation between exterior measurements, unique continuation, and
Runge approximation is central to inverse problems for nonlocal wave
equations. The preceding article \cite{zimmermann-paperI} develops the
finite energy theory for unequal elastic and dissipative orders and
uses approximation in $L^2(0,T;\widetilde H^r(\Omega))$ to recover a
linear potential. The present work addresses stronger approximation,
including the time derivative. Its smallness condition concerns the
rough-source linear theory; it is not needed in the finite energy
well-posedness result recalled from that article.

For coincident powers, the viscous wave equation in
\cite{zimmermann2024-viscous-wave} provides a related setting for
spatial-energy approximation and recovery of linear and nonlinear
perturbations. The work of Lin, Tyni, and Zimmermann
\cite{Semilinear-nonlocal-wave-eq} develops forward and inverse theory
for semilinear nonlocal wave equations. The article \cite{Optimal-Runge-nonlocal-wave} derives stronger Runge approximation property, which is then extended in \cite{zimmermann2024optimalrungeapproximationdamped} to damped nonlocal wave equations. Moreover, the article \cite{Third-order-nonlocal-wave} generalizes the work \cite{Semilinear-nonlocal-wave-eq} to third-order in time nonlocal wave equations with more general nonlinearities. These results motivate the use of nonlinear response identities together with Runge approximation. Here the nonlinearity acts on the velocity, so the recovery argument uses approximation in a space that controls time derivatives. The nonnegative damping structure also supplies a direct weighted monotonicity proof of the small-amplitude response.


\subsection{Organization of the article}
Section~\ref{R:sec: forward} recalls the finite energy theory, introduces
the space--time spaces, and describes the transposition identities.
Section~\ref{F:sec: very weak solutions} constructs very weak solutions,
first for spatial Sobolev sources and then for general energy-dual
sources. Section~\ref{R:sec: Runge for constant coeff} proves the
spatial-energy and full space--time Runge approximation theorems.
Section~\ref{sec:nonlinear-velocity-damping} develops the nonlinear
forward theory, proves the small-amplitude response formula, and
establishes uniqueness of the damping coefficient. Finally,
Appendix~\ref{R:sec: smoothing higher order heat eq} records a smoothing
estimate for higher-order heat equations.

\section{Forward theory and transposition identities}
\label{R:sec: forward}

In this section, we collect the forward well-posedness results and functional analytic tools that underlie the space--time Runge approximation argument developed later in the paper. In particular, we introduce the natural energy spaces for weakly viscous nonlocal waves, discuss weak and very weak solution concepts, and heuristically derive the transposition identities that form the basis of the duality argument.

\subsection{Fractional Laplacians and Bessel potential spaces}
\label{R:subsec: frac lap}

For $s\ge 0$, the fractional Laplacian $(-\Delta)^s$ is defined as the Fourier multiplier
\begin{equation}
\label{R:eq:fracLapFourDef}
(-\Delta)^s u
= \ifourier\bigl(|\xi|^{2s}\widehat{u}(\xi)\bigr),
\end{equation}
whenever the right-hand side is well defined. It is well known that
\[
(-\Delta)^s\colon H^{r,p}(\R^n)\to H^{r-2s,p}(\R^n)
\]
is bounded for all $r\in\R$ and $1<p<\infty$. This follows from the Mikhlin multiplier theorem applied to the symbol
\[
a(\xi)=\frac{|\xi|^{2s}}{\langle\xi\rangle^{2s}}.
\] 
We recall that the Bessel potential spaces $H^{r,p}(\R^n)$, $r\in\R$, $1\leq p<\infty$, are given by
\[
H^{r,p}(\R^n)
=
\{u\in \tempered(\R^n)\,;\,
\|u\|_{H^{r,p}(\R^n)}
\vcentcolon=
\|\langle D\rangle^r u\|_{L^p(\R^n)}
<\infty
\}.
\]
Moreover, for any open set $V\subset\R^n$, we define the closed subspace
\begin{equation}
\label{R:eq: local bessel}
    \widetilde{H}^{s,p}(V)\vcentcolon=
\overline{C_c^\infty(V)}^{\,H^{s,p}(\mathbb{R}^n)}.
\end{equation}
If $p=2$, then we shall denote the above spaces as $H^s(\R^n)$ and $\widetilde{H}^s(V)$.

\subsection{Weak solutions}
\label{R:subsec: weak sols}

In the remainder of this section, we discuss weak and very weak solutions for the weakly viscous nonlocal wave equation
\begin{equation}
\label{R:eq: weakly viscous wave forward}
\begin{cases}
\big( \partial_t^2 + (-\Delta)^r \partial_t + (-\Delta)^s +q\big) u = f & \text{in } \Omega_T, \\
u = 0 & \text{in } (\Omega_e)_T, \\
u(0) = u_0, \ \partial_t u(0) = u_1 & \text{in } \Omega.
\end{cases}
\end{equation}
Here and throughout this section, $\Omega\subset\R^n$ always denotes a bounded Lipschitz domain, $\Omega_e=\R^n\setminus\overline{\Omega}$ is the exterior of $\Omega$, $T>0$ is a finite time-horizon, and $0\leq r<s<\infty$ are two fractional exponents.

The next result introduces the class of potentials we will consider but also establishes its continuity properties.

 \begin{lemma}[{\cite[Theorem~3.1]{zimmermann2024-viscous-wave}}]
    \label{R:lemma: cond on q}
    Let $\mathcal{Q}_s$ be the class of all real-valued functions $q\in L^1_{loc}(\Omega_T)$ satisfying the following conditions:
    \begin{enumerate}[(i)]
        \item\label{R:item: integrability q} $q\in L^{\infty}(0,T;L^p(\Omega))$ for some exponent $2\leq p\leq \infty$ obeying
    \begin{equation}
    \label{R:eq: exponent p}
        \begin{cases}
			n/s\leq p\leq \infty, &\, \text{if }\, 2s< n,\\
			2<p\leq \infty,  &\, \text{if }\, 2s= n,\\
			2\leq p\leq \infty, &\, \text{if }\, 2s > n.
		\end{cases}
    \end{equation}
        \item $q$ is \emph{distributionally continuous in $t$}, meaning that
    \begin{equation}
    \label{R:eq: continuity in t}
        [0,T]\ni t\mapsto \int_{\Omega}q(x,t)\varphi(x)\,dx\in C([0,T])
    \end{equation}
    for all $\varphi\in C_c^{\infty}(\Omega)$.
    \end{enumerate}
    We denote by $\mathcal Q_s^p$ the subclass of all $q\in\mathcal Q_s$ satisfying condition~\ref{R:item: integrability q} with exponent $p$. 
    
    For each $q\in \mathcal{Q}_s^p$ and all $u,v\in \widetilde{H}^s(\Omega)$, we have
\[
t\mapsto \langle q(t)u,v\rangle_{L^2(\Omega)}\in C([0,T])
\]
and
\begin{equation}
\label{R:eq: estimate a1}
\left|
\langle q(t)u,v\rangle_{L^2(\Omega)}
\right|
\lesssim
\|q(t)\|_{L^p(\Omega)}
\|u\|_{\widetilde H^s(\Omega)}
\|v\|_{L^2(\Omega)}
\end{equation}
for a.e. $0\leq t\leq T$. The implicit constant depends only on $n,s,p$, and $\Omega$.
    \end{lemma}

\begin{theorem}[Weak solutions]
\label{R:thm: weak sol}
   Assume that $q\in\mathcal Q_s^p$. Suppose that the initial data $(u_0,u_1)$, and source function $f$ satisfy
   \begin{equation}
   \label{R:eq: data weak sol}
       (u_0,u_1,f)\in \widetilde{H}^s(\Omega)\times \widetilde{L}^2(\Omega)\times  L^2(0,T;H^{-r}(\Omega)).
   \end{equation}
   Then problem \eqref{R:eq: weakly viscous wave forward}
admits a unique finite energy solution $u\colon \R^n\times (0,T)\to\C$. More concretely, 
   \begin{enumerate}[(i)]
       \item\label{R:regularity} $u$ belongs to the space
       \[
     C([0,T];\widetilde H^s(\Omega)) \cap C^1([0,T];\widetilde L^2(\Omega)) \cap H^1(0,T;\widetilde H^r(\Omega)),
    \]
       \item\label{R:weak sol} $u$ weakly solves \eqref{R:eq: weakly viscous wave forward}, meaning that,
\begin{equation}
\label{R:eq: weak formulation}
\begin{split}
&\int_0^T
\left(
-\langle \partial_t u,\partial_t\varphi\rangle_{L^2(\Omega)}
+
\langle (-\Delta)^{r/2} \partial_t u,(-\Delta)^{r/2}\varphi\rangle_{L^2(\R^n)}
\right)\,dt
\\
&\quad
+
\int_0^T
\left(
\langle (-\Delta)^{s/2}u,(-\Delta)^{s/2}\varphi\rangle_{L^2(\R^n)}
+
\langle qu,\varphi\rangle_{L^2(\Omega)}
\right)\,dt
\\
&=
\langle u_1,\varphi(0)\rangle_{L^2(\Omega)}
+
\int_0^T
\langle f,\varphi\rangle_{H^{-r} (\Omega)\times\widetilde H^r(\Omega)}\,dt
\end{split}
\end{equation}
       for any $\varphi\in C_c^{\infty}(\Omega\times [0,T))$,
       \item\label{R:initial conds} and the initial conditions are attained pointwise
       \[
(u(0),\partial_t u(0))
=
(u_0,u_1)
\text{ in }
\widetilde H^s(\Omega)\times\widetilde L^2(\Omega).
\]
   \end{enumerate}
\end{theorem}

\begin{remark}
    If condition~\ref{R:regularity} is replaced by
    \[
       u\in L^{2}(0,T;\widetilde{H}^s(\Omega))\cap H^1(0,T;\widetilde{H}^r(\Omega))
    \]
    and we only assume that $u(0)=u_0$ in $\widetilde{H}^r(\Omega)$ instead of \ref{R:initial conds}, then we say $u$ is a \emph{weak finite energy solution} of \eqref{R:eq: weakly viscous wave forward}.
\end{remark}

\begin{remark}[{Non-homogeneous problems, \cite[Theorem~3.24]{zimmermann-paperI}}]
\label{R:Remark: non-homogeneous problems}
        Suppose that
         \[
            (u_0,u_1,\varphi)\in H^s(\R^n)\times L^2(\R^n)\times C^2([0,T];H^{2s}(\R^n))
         \]
        satisfies the compatibility condition 
        \begin{equation}
        \label{R:eq: compatibility cond}
            \varphi(0)=u_0\text{ and }\partial_t\varphi(0)=u_1\text{ a.e. in }\Omega_e.
        \end{equation}
        Then problem 
        \begin{equation}
\label{R:eq: weakly viscous wave non-hom forward}
\begin{cases}
\big( \partial_t^2 + (-\Delta)^r \partial_t + (-\Delta)^s +q\big) u = 0 & \text{in } \Omega_T, \\
u = \varphi & \text{in } (\Omega_e)_T, \\
u(0) = u_0, \ \partial_t u(0) = u_1 & \text{in } \Omega.
\end{cases}
\end{equation}
        has a unique finite energy solution $u$, that is,
        \begin{enumerate}[(a)]
            \item $u\in C([0,T];H^s(\R^n))\cap C^1([0,T];L^2(\R^n))\cap H^1(0,T;H^r(\R^n))$,
            \item  $u$ solves \eqref{R:eq: weakly viscous wave non-hom forward} in the weak sense; see equation~\eqref{R:eq: weak formulation}, 
            \item $u=\varphi$ in $(\Omega_e)_T$,
            \item and $(u,\partial_t u)(0)=(u_0,u_1)$ in $H^s(\R^n)\times L^2(\R^n)$.
        \end{enumerate}
\end{remark}

\subsection{Energy spaces}
\label{R:subsec: energy spaces}

We now introduce the space--time energy spaces that will be used to formulate very weak solutions.

\begin{definition}[Energy spaces]
\label{R:def: energy spaces}
We define
\begin{equation}
\label{R:eq: intersection space}
X^1_{s,r}(\Omega_T)
\vcentcolon=
L^2(0,T;\widetilde{H}^s(\Omega))
\cap
H^1(0,T;\widetilde{H}^r(\Omega)),
\end{equation}
endowed with the Hilbert norm induced by the inner product
\begin{equation}
\label{R:eq: intersection space norm}
\langle u,v\rangle_{X^1_{s,r}(\Omega_T)}
\vcentcolon=
\langle u,v\rangle_{L^2(0,T;\widetilde{H}^s(\Omega))}
+
\langle \partial_t u,\partial_tv\rangle_{L^2(0,T;\widetilde{H}^r(\Omega))}.
\end{equation}
We also introduce the terminal-value subspace
\begin{equation}
\label{R:eq: subspace}
X^1_{s,r,T}(\Omega_T)
\vcentcolon=
L^2(0,T;\widetilde{H}^s(\Omega))
\cap
H^1_T(0,T;\widetilde{H}^r(\Omega)),
\end{equation}
where $H^1_T(0,T;\widetilde{H}^r(\Omega))$ is the closed subspace of functions vanishing at time $T$, that is
\begin{equation}
\label{R:eq: vanishing at T}
H^1_T(0,T;\widetilde{H}^r(\Omega))
\vcentcolon=
\{u\in H^1(0,T;\widetilde{H}^r(\Omega))\,;\,u(T)=0\}.
\end{equation}
\end{definition}

The space $X^1_{s,r}(\Omega_T)$ is the natural energy space for weakly viscous equations in the regime $r < s$. Indeed, as seen in the previous section, weak finite energy solutions typically satisfy
\[
u\in L^2(0,T;\widetilde{H}^s(\Omega)),
\quad
\partial_t u\in L^2(0,T;\widetilde{H}^r(\Omega)).
\]
The subspace $X^1_{s,r,T}(\Omega_T)$ will be used for adjoint problems and transposition formulations, where terminal conditions are imposed at time $T$.

The next lemma collects the basic structural properties of these spaces, which is established in \cite[Lemma~2.5]{zimmermann-paperI}. We remark that throughout this work all Hilbert and Banach spaces are assumed to be complex. In particular, scalar products are taken to be linear in the first argument and antilinear in the second. Furthermore, if $X$ is a complex Banach space, then $X'$ denotes its antidual space, that is, the space of all continuous antilinear functionals on $X$.

\begin{lemma}[Properties of the energy spaces]
\label{R:lemma: space of sources}
The spaces $X^1_{s,r}(\Omega_T)$ and $X^1_{s,r,T}(\Omega_T)$ have the following properties.

\begin{enumerate}[(i)]
\item
\label{R:Banach prop of X}
Both $X^1_{s,r}(\Omega_T)$ and $X^1_{s,r,T}(\Omega_T)$ are Hilbert spaces. Moreover, the norm of $X^1_{s,r}(\Omega_T)$ is equivalent to the standard intersection norm
\[
\max\big(\|u\|_{L^2(0,T;\widetilde{H}^s(\Omega))},\|u\|_{H^1(0,T;\widetilde{H}^r(\Omega))}\big).
\]

\item
\label{R:Dual of X}
The antidual spaces can be identified as
\begin{equation}
\label{R:eq: dual space}
(X^1_{s,r}(\Omega_T))'
=
L^2(0,T;H^{-s}(\Omega))
+
\bigl(H^1(0,T;\widetilde{H}^r(\Omega))\bigr)'
\end{equation}
and
\begin{equation}
\label{R:eq: dual space X T}
(X^1_{s,r,T}(\Omega_T))'
=
L^2(0,T;H^{-s}(\Omega))
+
\bigl(H^1_T(0,T;\widetilde{H}^r(\Omega))\bigr)'.
\end{equation}
The sums are understood as sums of compatible Banach spaces and are endowed with the norms
\begin{equation}
\label{R:eq: norm dual X}
\inf_{F=F_1+F_2}
\left(
\|F_1\|_{L^2(0,T;H^{-s}(\Omega))}
+
\|F_2\|_{\bigl(H^1(0,T;\widetilde{H}^r(\Omega))\bigr)'}
\right)\sim \|F\|_{(X^1_{s,r}(\Omega_T))'}
\end{equation}
and
\begin{equation}
\label{R:eq: norm dual X T}
\inf_{F=F_1+F_2}
\left(
\|F_1\|_{L^2(0,T;H^{-s}(\Omega))}
+
\|F_2\|_{\bigl(H^1_T(0,T;\widetilde{H}^r(\Omega))\bigr)'}
\right)\sim \|F\|_{(X^1_{s,r,T}(\Omega_T))'}.
\end{equation}

\item
\label{R:Density in dual of X}
The spaces
\[
X^1_{s,r}(\Omega_T),
\qquad
L^2(0,T;\widetilde{L}^2(\Omega)),
\qquad
(X^1_{s,r}(\Omega_T))'
\]
form a Gelfand triple, where the second embedding is induced by the $L^2$ inner product and the identification of a Hilbert space with its antidual via the Riesz representation theorem. In other words, the continuous embeddings
\begin{equation}
\label{R:eq: dense embedding}
X^1_{s,r}(\Omega_T)
\hookrightarrow
L^2(0,T;\widetilde{L}^2(\Omega))
\hookrightarrow
(X^1_{s,r}(\Omega_T))'
\end{equation}
are dense.

\item
\label{R:density in X T}
Let $(\varphi_k)_{k\in\N}\subset \widetilde{H}^s(\Omega)$ be the eigenfunctions of the fractional Laplacian $(-\Delta)^s$ with homogeneous exterior condition. Then the set
\begin{equation}
\label{R:eq: dense set D T}
\mathcal{D}_T
\vcentcolon=
\left\{
\sum_{k=1}^N \eta_k(t)\varphi_k(x)
\,;\,
\eta_k\in C_c^\infty([0,T)),\
1\leq k\leq N,\
N\in\N
\right\}
\end{equation}
is dense in $X^1_{s,r,T}(\Omega_T)$.
\end{enumerate}
\end{lemma}

\begin{remark}[Initial-value subspace]
    Similarly, we may define initial-value subspace
\begin{equation}
\label{R:eq: energy subspace vanishing 0}
    X^1_{s,r,0}(\Omega_T)=L^2(0,T;\widetilde{H}^s(\Omega))\cap H^1_0(0,T;\widetilde{H}^r(\Omega)),
\end{equation}
where $H^1_0(0,T;\widetilde{H}^r(\Omega))$ is the closed subspace of functions vanishing at time $0$, that is
\begin{equation}
\label{R:eq: vanishing at 0}
    H^1_0(0,T;\widetilde{H}^r(\Omega))
\vcentcolon=
\{u\in H^1(0,T;\widetilde{H}^r(\Omega))\,;\,u(0)=0\}.
\end{equation}
    All properties of $X^1_{s,r,T}(\Omega_T)$ carry over to $X^1_{s,r,0}(\Omega_T)$, while the only difference being that we take in \eqref{R:eq: dense set D T} the functions $\eta_k$ in  $C_c^{\infty}((0,T])$ instead of $C^{\infty}_c([0,T))$.
\end{remark}

\subsection{Very weak solutions}
\label{R:subsec: very weak}

Now we can formulate the well-posedness result for weakly viscous equations in the class of very weak solutions. That is, we want to introduce solutions $u$ to 
\begin{equation}
\label{R:eq: very weakly viscous wave forward}
\begin{cases}
\big( \partial_t^2 + (-\Delta)^r \partial_t + (-\Delta)^s +q\big) u = f & \text{in } \Omega_T, \\
u = 0 & \text{in } (\Omega_e)_T, \\
u(0) = u_0, \ \partial_t u(0) = u_1 & \text{in } \Omega,
\end{cases}
\end{equation}
when the initial data $(u_0,u_1)$ and sources $f$ only satisfy
    \begin{equation}
    \label{R:eq: rough data}
(u_0,u_1)\in
\begin{cases}
\widetilde H^{2r-s}(\Omega)\times H^{-s}(\Omega),
& s<2r,\\[2mm]
\widetilde L^2(\Omega)\times H^{-s}(\Omega),
& 2r\leq s,
\end{cases}
\qquad f\in (X^1_{s,r}(\Omega_T))'.
\end{equation}
In view of Section~\ref{R:subsec: weak sols}, we observe that the initial data and sources are too irregular to allow the existence of weak solutions.

The construction of very weak solutions is based on a transposition argument using the adjoint terminal-value problem, an approximation of the rough data $(u_0,u_1,f)$ by regular data
\[
    (u_{0,n},u_{1,n},f_n)\in \widetilde{H}^s(\Omega)\times\widetilde{L}^2(\Omega)\times L^2(0,T;H^{-r}(\Omega)),
\]
and the existence of finite energy solutions $u_n$ to \eqref{R:eq: very weakly viscous wave forward} for the data together with suitable energy estimates.

\begin{theorem}[Very weak solutions]
\label{R:thm: very weak sol}
    There exists a constant
$C_\ast>0$, depending only on $n,r,s,p$, and $\Omega$, such that if
$q\in \mathcal{Q}_s^p$ satisfies
\begin{equation}
\label{R:eq: smallness condition on q rough}
    \|q\|_{L^1(0,T;L^p(\Omega))}
    \leq \frac{1}{8C_\ast},
\end{equation}
then the following assertions hold.

For all data $(u_0,u_1,f)$ satisfying \eqref{R:eq: rough data}, there exists a unique very weak solution
\begin{equation}
\label{R:eq: reg very weak sol}
u\in
L^\infty(0,T;\widetilde L^2(\Omega))
\cap
L^2(0,T;\widetilde H^r(\Omega))
\end{equation}
of \eqref{R:eq: very weakly viscous wave forward}, in the sense that for every $g\in L^2(0,T;H^{-r}(\Omega))$, with $v\in  X^1_{s,r,T}(\Omega_T)$ denoting the finite energy
solution of the adjoint terminal-value problem 
\begin{equation}
\label{R:eq: adjoint problem very weak}
\begin{cases}
(\partial_t^2-(-\Delta)^r \partial_t+(-\Delta)^s+q)v=g
&\text{in }\Omega_T,
\\
v=0
&\text{in }(\Omega_e)_T,
\\
v(T)=0,\ \partial_tv(T)=0
&\text{in }\Omega,
\end{cases}
\end{equation}
one has the transposition identity
\begin{equation}
\label{R:eq: identity for very weak sol rough}
\begin{split}
&\langle f,v\rangle_{(X^1_{s,r}(\Omega_T))'\times X^1_{s,r}(\Omega_T)}
=
\int_0^T\langle u,g\rangle_{H^{-r}(\Omega)\times \widetilde H^r(\Omega)}\,dt
+
\langle u_0,\partial_t v(0)\rangle_{L^2(\Omega)}\\
&\quad
-
\langle u_1,v(0)\rangle_{H^{-s}(\Omega)\times\widetilde H^s(\Omega)}-
\langle
(-\Delta)^{(2r-s)/2}u_0,
(-\Delta)^{s/2}v(0)
\rangle_{L^2(\R^n)}
\end{split}
\end{equation}
when $s<2r$, and
\begin{equation}
\label{R:eq: identity for very weak sol rough 2}
\begin{split}
&\langle f,v\rangle_{(X^1_{s,r}(\Omega_T))'\times X^1_{s,r}(\Omega_T)}
=
\int_0^T\langle u,g\rangle_{H^{-r}(\Omega)\times \widetilde H^r(\Omega)}\,dt
+
\langle u_0,\partial_t v(0)\rangle_{L^2(\Omega)}
\\
&\quad
-
\langle u_1,v(0)\rangle_{H^{-s}(\Omega)\times\widetilde H^s(\Omega)}-
\langle u_0,(-\Delta)^r v(0)\rangle_{L^2(\Omega)}
\end{split}
\end{equation}
when $2r\leq s$.
\end{theorem}

\begin{remark}[Time-reversal and adjoint equation]
\label{R:remark: time reversal and adjoint equation}
Let $X$ be a Banach space and let $h\colon (0,T)\to X$. We define the \emph{time-reversal} of $h$ by
\begin{equation}
\label{R:eq: time reversal}
    h^\star(t)\vcentcolon = h(T-t),
    \qquad 0<t<T.
\end{equation}
In particular, the time-reversal operator interchanges the energy spaces with vanishing initial and terminal conditions, that is,
\begin{equation}
\label{R:eq: time reversal of energy space}
    (X^1_{s,r,T}(\Omega_T))^\star
    =
    X^1_{s,r,0}(\Omega_T),
    \qquad
    (X^1_{s,r,0}(\Omega_T))^\star
    =
    X^1_{s,r,T}(\Omega_T).
\end{equation}

Moreover, the unique finite energy solution of
\eqref{R:eq: adjoint problem very weak}
is given by
\[
    v(x,t)=\omega^\star(x,t),
\]
where $\omega$ denotes the unique finite energy solution of
\begin{equation}
\label{R:eq: transformed adjoint problem}
\begin{cases}
(\partial_t^2+(-\Delta)^r\partial_t+(-\Delta)^s+q^\star)\omega
=
g^\star
&\text{in }\Omega_T,
\\
\omega=0
&\text{in }(\Omega_e)_T,
\\
\omega(0)=0,\ \partial_t\omega(0)=0
&\text{in }\Omega,
\end{cases}
\end{equation}
whose existence and uniqueness follow from
Theorem~\ref{R:thm: weak sol}.
\end{remark}

\subsection{Transposition identities and duality pairing}
\label{R:sec: transposition}

   For the reader's convenience, we explain here the formal duality argument underlying the definition of very weak solutions.
   
   First, we observe that identities
\eqref{R:eq: identity for very weak sol rough}
and
\eqref{R:eq: identity for very weak sol rough 2}
are obtained formally by multiplying
\eqref{R:eq: very weakly viscous wave forward}
with the adjoint solution $v$ and integrating by parts in both space and time. Indeed, suppose for the sake of motivation that all functions are real-valued and sufficiently regular, then 
\[
\begin{split}
    \int_{\Omega_T}f v\,dxdt&=\int_{\Omega_T}[(\partial_t^2+(-\Delta)^r \partial_t+(-\Delta)^s+q)u]v\,dxdt\\
    &=\int_{\Omega_T}u(\partial_t^2-(-\Delta)^r \partial_t+(-\Delta)^s+q)v\,dxdt\\
    &\quad +\left.\int_{\Omega}\partial_t u v\,dx\right|_{t=0}^{t=T}-\left.\int_{\Omega}u\partial_t v\,dx\right|_{t=0}^{t=T}+\left.\int_{\Omega} [(-\Delta)^r u]v\,dx\right|_{t=0}^{t=T}\\
    &=\int_{\Omega_T}u g\,dxdt-\int_{\Omega}u_1 v(0)\,dx+\int_{\Omega}u_0\partial_t v(0)\,dx-\int_{\Omega} [(-\Delta)^r u_0]v(0)\,dx.
\end{split}
\]
In the first equality we used the equation for $u$, in the second we performed an integration by parts and in the last equality we used the equation for $v$ together with initial and final conditions of $u$ and $v$, respectively. In conclusion, we have the integral identity
\begin{equation}
\label{R:eq: integral identity}
\int_{\Omega_T}f v\,dxdt=\int_{\Omega_T}u g\,dxdt-\int_{\Omega}u_1 v(0)\,dx+\int_{\Omega}u_0\partial_t v(0)\,dx-\int_{\Omega} [(-\Delta)^r u_0]v(0)\,dx.
\end{equation}
On the level of smooth solutions and data, this already justifies the transposition identities but in the setup of Theorem~\ref{R:thm: very weak sol} the appearing functions only satisfy:
\begin{enumerate}[(i)]
    \item\label{R:item 1} $u\in L^\infty(0,T;\widetilde L^2(\Omega))
\cap
L^2(0,T;\widetilde H^r(\Omega))$,
    \item\label{R:item 2} $(u_0,u_1,f)$ satisfy
    \[
        (u_0,u_1)\in
\begin{cases}
\widetilde H^{2r-s}(\Omega)\times H^{-s}(\Omega),
& s<2r,\\[2mm]
\widetilde L^2(\Omega)\times H^{-s}(\Omega),
& 2r\leq s,
\end{cases}
\qquad f\in (X^1_{s,r}(\Omega_T))',
    \]
    \item\label{R:item 3} $v\in  C([0,T];\widetilde H^s(\Omega)) \cap C^1([0,T];\widetilde L^2(\Omega)) \cap H^1(0,T;\widetilde H^r(\Omega)) \subset X^1_{s,r,T}(\Omega_T)$,
    \item\label{R:item 4} $(v(0),\partial_t v(0),g)\in \widetilde{H}^s(\Omega)\times \widetilde{L}^2(\Omega)\times L^2(0,T;H^{-r}(\Omega))$.
\end{enumerate}
Hence, the term 
\[
    \int_{\Omega} [(-\Delta)^r u_0]v(0)\,dx
\]
is not necessarily well-defined, thus we write it, for $2r\leq s$, in the form 
\[
\int_{\Omega}u_0(-\Delta)^r v(0)\,dx
\]
and, for $2r>s$, as
\[
    \int_{\R^n}[(-\Delta)^{(2r-s)/2}u_0][
(-\Delta)^{s/2}v(0)]\,dx.
\]
All remaining terms in \eqref{R:eq: integral identity} are well defined by properties \ref{R:item 1}--\ref{R:item 4} together with the Gelfand triples
\[
\widetilde H^r(\Omega)\hookrightarrow \widetilde L^2(\Omega)\hookrightarrow H^{-r}(\Omega)
\]
and
\[
X^1_{s,r}(\Omega_T)
\hookrightarrow
L^2(0,T;\widetilde L^2(\Omega))
\hookrightarrow
(X^1_{s,r}(\Omega_T))',
\]
whose embeddings are continuous and dense; see equation~\eqref{R:eq: local bessel} and Lemma~\ref{R:lemma: space of sources} \ref{R:Density in dual of X}.

\section{Very weak theory}
\label{F:sec: very weak solutions}

In Sections~\ref{F:sec: very weak constant coeff}--\ref{F:sec: rough sources constant}, we construct very weak solutions to weakly viscous
nonlocal wave equations
\begin{equation}
\label{F:eq: very weak problem constant coeff}
\begin{cases}
(\partial_t^2+(-\Delta)^r \partial_t+(-\Delta)^s+q)u=f
&\text{in }\Omega_T,
\\
u=0
&\text{in }(\Omega_e)_T,
\\
u(0)=u_0,
\quad
\partial_t u(0)=u_1
&\text{in }\Omega
\end{cases}
\end{equation}
for rough source terms $f$ and low regularity
initial data $(u_0,u_1)$. 

Our main objective is the development of a transposition framework for source terms belonging to the antidual space $(X^1_{s,r,T}(\Omega_T))'$.
To simplify several technical aspects of the construction, we first
study $f\in L^2(0,T;H^{-s}(\Omega))$,
which serves as a preparatory model case (see Section~\ref{F:sec: very weak constant coeff}). The corresponding estimates and approximation arguments will later allow us to streamline the
construction of very weak solutions in the more general transposition
framework. The results developed in this section require that the potential is small. 

The main difficulty compared to the classical nonlocal wave equation
is that the viscous term $(-\Delta)^r \partial_t u$
cannot be treated as a lower-order perturbation of $\partial_t^2+(-\Delta)^s$
in the weakly viscous regime $0<r<s$. Our construction therefore uses the parabolic regularization from
\cite[Section~3.1]{zimmermann-paperI} together with a double approximation
argument: first the rough data are approximated by finite energy data, and then
the resulting finite energy solutions are obtained as limits of their parabolically
regularized counterparts.

In
\cite{zimmermann-paperIII}, we generalize the theory to equations with variable coefficients.

Finally, the article \cite{zimmermann-paperIV} develops an alternative theory of very weak solutions of damped nonlocal wave equations with sources in  $L^2(0,T;H^{-s}(\Omega))$. The main advantage of this theory is the fact that it does not require any smallness assumption of the potential.

Throughout this section, $\Omega\subset\R^n$ will always denote a bounded Lipschitz domain, $T>0$, and
$0\leq r<s<\infty$, if not specified further.

\subsection{Very weak theory for equations with 
\texorpdfstring{$L^2(0,T;H^{-s}(\Omega))$}{L2 H-minus-s}-sources}
\label{F:sec: very weak constant coeff}

First let us note that for \eqref{F:eq: very weak problem constant coeff} with rough data
\begin{equation}
\label{F:eq: data H minus s sol}
    (u_0,u_1,f)\in \widetilde H^r(\Omega)\times H^{-s}(\Omega)\times L^2(0,T;H^{-s}(\Omega)),
\end{equation}
finite energy solutions need not exist. We therefore define
solutions by duality against solutions of the adjoint terminal-value problem
\begin{equation}
\label{F:eq: adjoint problem very weak}
\begin{cases}
(\partial_t^2-(-\Delta)^r \partial_t+(-\Delta)^s+q)v=g
&\text{in }\Omega_T,
\\
v=0
&\text{in }(\Omega_e)_T,
\\
v(T)=0,\ \partial_tv(T)=0
&\text{in }\Omega.
\end{cases}
\end{equation}
Equivalently, after the time reversal $w(t)=v(T-t)$, the potential is replaced by
$q^\star(x,t)=q(x,T-t)$.

By the finite energy theory of \cite[Theorem~3.21]{zimmermann-paperI}, for every
$g\in L^2(0,T;H^{-r}(\Omega))$, the adjoint problem
\eqref{F:eq: adjoint problem very weak} admits a unique finite energy solution
\[
v\in
C([0,T];\widetilde H^s(\Omega))
\cap
C^1([0,T];\widetilde L^2(\Omega))
\cap
H^1(0,T;\widetilde H^r(\Omega)).
\]
Multiplying equation \eqref{F:eq: very weak problem constant coeff} by $v$ and carrying out a formal integration by parts, motivates the following definition of very weak solutions for rough initial data and sources.

\begin{definition}[Very weak solutions for Sobolev sources]
\label{F:def: very weak solution L2 Hminus s}
Let $(u_0,u_1,f)$ satisfy \eqref{F:eq: data H minus s sol}.
A function
\[
u\in
L^\infty(0,T;\widetilde L^2(\Omega))
\cap
L^2(0,T;\widetilde H^r(\Omega))
\]
is called \emph{very weak solution} of \eqref{F:eq: very weak problem constant coeff} if
for every source $g\in L^2(0,T;H^{-r}(\Omega))$, with $v$ denoting the finite energy
solution of \eqref{F:eq: adjoint problem very weak}, one has
\begin{equation}
\label{F:eq: identity for very weak sol intro}
\begin{split}
&\int_0^T\langle f,v\rangle_{H^{-s}(\Omega)\times\widetilde H^s(\Omega)}\,dt
=
\int_0^T\langle u,g\rangle_{H^{-r}(\Omega)\times\widetilde H^r(\Omega)}\,dt
+
\langle u_0,\partial_t v(0)\rangle_{L^2(\Omega)}
\\
&\quad
-
\langle u_1,v(0)\rangle_{H^{-s}(\Omega)\times\widetilde H^s(\Omega)}-
\langle
(-\Delta)^{r/2}u_0,
(-\Delta)^{r/2}v(0)
\rangle_{L^2(\R^n)}.
\end{split}
\end{equation}
\end{definition}

Our well-posedness result reads as follows.

\begin{theorem}[Very weak solutions with sources in
\texorpdfstring{$L^2(0,T;H^{-s}(\Omega))$}{L2 H-minus-s}]
\label{F:thm: ver weak sol}
    There exists a constant
$C_\ast>0$, depending only on $n,r,s,p$, and $\Omega$, such that if
$q\in \mathcal{Q}_s^p$ satisfies
\begin{equation}
\label{F:eq: smallness condition on q and f}
    \|q\|_{L^1(0,T;L^p(\Omega))}
    \leq
    \frac{1}{8C_\ast},
\end{equation}
then, for all $(u_0,u_1,f)$ obeying \eqref{F:eq: data H minus s sol},
there exists a unique very weak solution
\[
u\in
L^\infty(0,T;\widetilde L^2(\Omega))
\cap
L^2(0,T;\widetilde H^r(\Omega))
\]
of \eqref{F:eq: very weak problem constant coeff} in the sense of Definition~\ref{F:def: very weak solution L2 Hminus s}. Moreover, $u$ satisfies
\begin{equation}
\label{F:eq: energy estimate very weak solutions}
\begin{split}
&
\|u(t)\|_{L^2(\Omega)}^2
+
\int_0^t
\|(-\Delta)^{r/2}u(\tau)\|_{L^2(\R^n)}^2\,d\tau
\\
&\lesssim
\left(
\|u_0\|_{\widetilde H^r(\Omega)}^2
+
\|u_1\|_{H^{-s}(\Omega)}^2
+
\|f\|_{L^2(0,t;H^{-s}(\Omega))}^2
\right)
\exp\left(
C\|q\|_{L^1(0,t;L^p(\Omega))}
\right)
\end{split}
\end{equation}
for a.e. $0<t<T$.
\end{theorem}

    \begin{proof}
    \emph{Step 1: Approximation and regularization.}
        By density, we may choose
        \[
            (u_{0,n},u_{1,n},f_n)\in \widetilde{H}^s(\Omega)\times\widetilde{L}^2(\Omega)\times L^2(0,T;H^{-r}(\Omega)),\quad n\in\N,
        \]
        such that 
        \[
         (u_{0,n},u_{1,n},f_n)\to (u_0,u_1,f) \text{ in }\widetilde{H}^r(\Omega)\times H^{-s}(\Omega)\times L^2(0,T;H^{-s}(\Omega))
        \]
        as $n\to\infty$. Next, let us consider the approximate problem
        \begin{equation}
        \label{F:eq: very weak sol approximated}
            \begin{cases}
				( \partial_t^2 +(-\Delta)^r \partial_t +(-\Delta)^s+q) u=f_n &\text{ in }\Omega_T,\\
				u=0     &\text{ in }(\Omega_e)_T,\\
				u(0)=u_{0,n},\ \partial_t u(0)=u_{1,n}  &\text{ in }\Omega.
			\end{cases}
        \end{equation}
        By \cite[Theorem~3.21]{zimmermann-paperI} and \cite[Theorem~3.22]{zimmermann-paperI}, this problem has a unique finite energy solution
        \[
            u_n\in C([0,T];\widetilde{H}^s(\Omega))\cap C^1([0,T];\widetilde{L}^2(\Omega))\cap H^1(0,T;\widetilde{H}^r(\Omega)).
        \]
        Furthermore, let 
        \[
            v\in C([0,T];\widetilde{H}^s(\Omega))\cap C^1([0,T];\widetilde{L}^2(\Omega))\cap H^1(0,T;\widetilde{H}^r(\Omega))
        \]
        be the unique finite energy solution to the adjoint problem \eqref{F:eq: adjoint problem very weak}, where $g\in L^2(0,T;H^{-r}(\Omega))$. Moreover, from \cite[Corollary~3.13]{zimmermann-paperI}, we know that for any $\eps>0$ there exist 
        \[
            u_{\eps,n}, v_\eps\in \widetilde{W}(0,T;\widetilde{H}^s(\Omega))
        \]
        solving the related parabolically regularized problems, where we take the viscosity term $+\eps (-\Delta)^s\partial_t$ for \eqref{F:eq: very weak sol approximated} and $-\eps(-\Delta)^s\partial_t$ for \eqref{F:eq: adjoint problem very weak} as the second is a terminal value problem.
        We assert that for each $n\in\N$ the sequence $(u_{\eps,n})_{\eps>0}$ has the following convergence and boundedness properties,
        \begin{enumerate}[(I)]
            \item\label{F:convergence 1 very} $u_{\eps,n}\weakstar u_n$ in $L^{\infty}(0,T;\widetilde{H}^s(\Omega))$,
            \item\label{F:convergence 2 very} $\partial_t u_{\eps,n}\weakstar \partial_t u_n$ in $L^{\infty}(0,T;\widetilde{L}^2(\Omega))$,
            \item\label{F:convergence 3 very} $\partial_t u_{\eps,n}\weak \partial_t u_n$ in $L^{2}(0,T;\widetilde{H}^r(\Omega))$,
            \item\label{F:convergence 4 very} $u_{\eps,n}(t)\to u_n(t)$ in $\widetilde{H}^s(\Omega)$ and $\partial_t u_{\eps,n}(t)\to \partial_t u_n(t)$ in $\widetilde{L}^2(\Omega)$ for $0\leq t\leq T$,
            \item\label{F:boundedness very} $\eps^{1/2}\partial_t u_{\eps,n}$ is bounded in $L^2(0,T;\widetilde{H}^s(\Omega))$.
        \end{enumerate}
       	The same holds for the sequence $(v_\eps)_{\eps>0}$ and $v$. The proofs of the related assertions are identical after time reversal. Let us note that only assertion \ref{F:convergence 4 very} is not an immediate consequence of the proof of \cite[Proposition~3.16]{zimmermann-paperI}. Furthermore, one has
        \begin{equation}
        \label{F:eq: strong convergence of damping term}
            \partial_t u_{\eps,n}\to \partial_t u_n\text{ in }L^{2}(0,T;\widetilde{H}^r(\Omega))
        \end{equation}
        as $\eps\to 0$. We postpone the proof of \ref{F:convergence 4 very} and \eqref{F:eq: strong convergence of damping term} until after the theorem. Moreover, we notice that the convergence \eqref{F:eq: strong convergence of damping term} yields
        \begin{equation}
        \label{F:eq: strong conv of w eps}
            u_{\eps,n}\to u_n\quad\text{in}\quad C^{0,1/2}([0,T];\widetilde{H}^r(\Omega))
        \end{equation}
        as $\eps\to 0$. 

        In fact, the fundamental theorem of calculus implies
        \[
            \|u_{\eps,n}-u_{n}\|_{L^{\infty}(0,T;\widetilde{H}^r(\Omega))}\leq T^{1/2}\|\partial_t u_{\eps,n} -\partial_t u_{n}\|_{L^2(0,T;\widetilde{H}^r(\Omega))}\to 0
        \]
       and
        \[
            [u_{\eps,n}-u_{n}]_{C^{0,1/2}([0,T];\widetilde{H}^r(\Omega))}\leq \|(-\Delta)^{r/2} \partial_t (u_{\eps,n} -u_{n})\|_{L^2(0,T;L^2(\R^n))}\to 0
        \]
        as $\eps\to 0$.\\
        
        \emph{Step 2: The transposition identity for smooth data.} Now, testing the PDE \eqref{F:eq: very weak sol approximated} for $u_{\eps,n}$ by $v_\eps$, we deduce
        \begin{equation}
        \label{F:eq: identity for regularized problem 1}
            \begin{split}
                \int_0^T\langle f_n,v_\eps\rangle\,dt&= \int_0^T\left(\langle \partial_t^2 u_{\eps,n},v_\eps\rangle+\langle (-\Delta)^{r/2} \partial_t u_{\eps,n} ,(-\Delta)^{r/2}v_\eps\rangle_{L^2(\R^n)}\right)\,dt\\
                &\quad +\int_0^T\left(\langle (-\Delta)^{s/2}u_{\eps,n},(-\Delta)^{s/2}v_\eps\rangle_{L^2(\R^n)}+\langle qu_{\eps,n},v_\eps\rangle_{L^2(\Omega)}\right)\,dt\\
                &\quad+\eps \int_0^T\langle (-\Delta)^{s/2}\partial_t u_{\eps,n},(-\Delta)^{s/2} v_\eps\rangle_{L^2(\R^n)}\,dt.
            \end{split}
        \end{equation}
        By the regularity properties of $u_{\eps,n}$ and $v_\eps$, we can apply the integration by parts formula, which yields
        \begin{equation}
        \label{F:eq: integration by parts second derivative}
            \begin{split}
                \int_0^T\langle \partial_t^2 u_{\eps,n},v_\eps\rangle\,dt&=\langle u_{0,n},\partial_t v_{\eps}(0)\rangle_{L^2(\Omega)}-\langle u_{1,n},v_\eps (0)\rangle_{L^2(\Omega)}+                \int_0^T\langle u_{\eps,n},\partial_t^2 v_{\eps}\rangle\,dt
            \end{split}
        \end{equation}
        and
        \begin{equation}
        \label{F:eq: integration by parts first derivative}
            \begin{split}
               & \int_0^T\langle (-\Delta)^{\alpha/2}\partial_t u_{\eps,n} ,(-\Delta)^{\alpha/2}v_\eps\rangle_{L^2(\R^n)}\,dt\\
               &=-\langle (-\Delta)^{\alpha/2}u_{0,n},(-\Delta)^{\alpha/2}v_\eps(0)\rangle_{L^2(\R^n)}\\
                &\quad-\int_0^T\langle (-\Delta)^{\alpha/2}u_{\eps,n},(-\Delta)^{\alpha/2}\partial_t v_{\eps}\rangle_{L^2(\R^n)}\,dt
            \end{split}
        \end{equation}
        for $\alpha=r$ or $\alpha=s$.
        Inserting \eqref{F:eq: integration by parts second derivative} -- \eqref{F:eq: integration by parts first derivative} (with $\alpha=r,s$) into \eqref{F:eq: identity for regularized problem 1} and using the PDE for $v_\eps$, we deduce
        \begin{equation}
        \label{F:eq: identity for regularized problem 2}
            \begin{split}
                \int_0^T\langle f_n,v_\eps\rangle\,dt&= \int_0^T \langle u_{\eps,n},g\rangle\,dt+\langle u_{0,n},\partial_t v_{\eps}(0)\rangle_{L^2(\Omega)}-\langle u_{1,n},v_\eps (0)\rangle_{L^2(\Omega)}\\
                &\quad -\langle (-\Delta)^{r/2}u_{0,n},(-\Delta)^{r/2}v_\eps(0)\rangle_{L^2(\R^n)}\\
                &\quad -\eps\langle (-\Delta)^{s/2}u_{0,n},(-\Delta)^{s/2}v_\eps(0)\rangle_{L^2(\R^n)}.
            \end{split}
        \end{equation}
        Next, we want to pass to the limit $\eps\to 0$. Observe that \ref{F:boundedness very}, $v_\eps(T)=0$ in $\widetilde{H}^s(\Omega)$, and the fundamental theorem of calculus guarantee that $\eps^{1/2}(-\Delta)^{s/2}v_\eps(0)$ is uniformly bounded in $L^2(\R^n)$. In fact, since $v_\eps\in H^1(0,T;\widetilde{H}^s(\Omega))$ implies
        \[
            (-\Delta)^{s/2}v_\eps (0)=-\int_0^T (-\Delta)^{s/2}\partial_t v_\eps\,d\tau
        \]
        in $L^2(\R^n)$ and so
        \[
        \begin{split}
             \|\eps^{1/2}(-\Delta)^{s/2} v_\eps (0)\|_{L^2(\R^n)}&
             \leq \int_0^T \|\eps^{1/2}(-\Delta)^{s/2}\partial_t v_\eps\|_{L^2(\R^n)}\,d\tau\\
             &\leq T^{1/2}\|\eps^{1/2}\partial_t v_\eps\|_{L^2(0,T;\widetilde{H}^s(\Omega))}.
        \end{split}
        \]
        Since the last expression is uniformly bounded in $\eps$, we have
        \begin{equation}
        \label{F:eq: conv regularization term}
            \eps\langle (-\Delta)^{s/2}u_{0,n},(-\Delta)^{s/2}v_\eps(0)\rangle_{L^2(\R^n)}\to 0\text{ as }\eps\to 0.
        \end{equation}
        Using the convergence results \ref{F:convergence 1 very}, \ref{F:convergence 4 very}, and \eqref{F:eq: conv regularization term}, we deduce from \eqref{F:eq: identity for regularized problem 2} the identity 
        \begin{equation}
       \label{F:eq: very weak formulation for fn}
              \begin{split}
                \int_0^T\langle f_n,v\rangle\,dt&= \int_0^T \langle u_{n},g\rangle\,dt+\langle u_{0,n},\partial_t v(0)\rangle_{L^2(\Omega)}-\langle u_{1,n},v(0)\rangle\\
                &\quad -\langle (-\Delta)^{r/2}u_{0,n},(-\Delta)^{r/2}v(0)\rangle_{L^2(\R^n)}.
            \end{split}
        \end{equation}
        
        \emph{Step 3: Uniform estimates.} For passing to the limit $n\to\infty$, we need to derive suitable energy estimates. We proceed similarly as in the uniqueness proof \cite[Theorem~3.22]{zimmermann-paperI} and consider the test function 
        \begin{equation}
            \label{F:eq: test function very weak sol}
           \varphi_n(t)=\begin{cases}
               -\int_{t}^{t_0}u_n(\tau)\,d\tau,&\text{ for }0\leq t\leq t_0,\\
               0,&\text{ for } t_0\leq t\leq T,
           \end{cases}
        \end{equation}
        where $0<t_0<T$ is fixed. By the arguments in the proof of \cite[Theorem~3.22]{zimmermann-paperI}, we can test the PDE \eqref{F:eq: very weak sol approximated} with $\varphi_n$, which yields
        \begin{equation}
          \begin{split}
                &\int_0^{t_0}\left(\partial_t\frac{\|u_n\|_{L^2(\Omega)}^2-\|(-\Delta)^{s/2}\varphi_n\|_{L^2(\R^n)}^2}{2}+\|(-\Delta)^{r/2}u_n\|_{L^2(\R^n)}^2-\re\langle qu_n,\varphi_n\rangle_{L^2(\Omega)}\right)\,dt\\
                &=\re\left(-\langle (-\Delta)^{r/2}u_{0,n},(-\Delta)^{r/2}\varphi_n(0)\rangle_{L^2(\R^n)}-\langle u_{1,n},\varphi_n(0)\rangle_{L^2(\Omega)}-\int_0^T\langle f_n,\varphi_n\rangle_{L^2(\Omega)}\,dt\right);
            \end{split}
        \end{equation}
        see \cite[Eq.~(3.13), Proof of Theorem~3.22]{zimmermann-paperI}.
        Hence, we deduce that
        \begin{equation}
        \label{F:eq: useful identity very weak}
            \begin{split}
                &\|u_n(t_0)\|_{L^2(\Omega)}^2+\|(-\Delta)^{s/2}\varphi_n(0)\|_{L^2(\R^n)}^2+2\int_0^{t_0}\|(-\Delta)^{r/2}u_n\|_{L^2(\R^n)}^2\,dt\\
                &=\|u_{0,n}\|_{L^2(\Omega)}^2-2\re\langle (-\Delta)^{r/2}u_{0,n},(-\Delta)^{r/2}\varphi_n(0)\rangle_{L^2(\R^n)}-2\re\langle u_{1,n},\varphi_n(0)\rangle_{L^2(\Omega)}\\
                &+2\re\int_0^{t_0}\left(\langle qu_n,\varphi_n\rangle_{L^2(\Omega)}-\langle f_n,\varphi_n\rangle_{L^2(\Omega)}\right)\,dt.
            \end{split}
        \end{equation}
        For any $\rho>0$, we can estimate
        \begin{equation}
        \label{F:eq: estimate source term}
        \begin{split}
            \left|\int_0^{t_0}\langle f_n,\varphi_n\rangle_{L^2(\Omega)}\,dt\right|&\leq \|f_n\|_{L^2(0,t_0;H^{-s}(\Omega))}\|\varphi_n\|_{L^2(0,t_0;\widetilde{H}^s(\Omega))}\\
            &\leq C_\rho \|f_n\|_{L^2(0,t_0;H^{-s}(\Omega))}^2+\rho\|(-\Delta)^{s/2}\varphi_n\|_{L^2(0,t_0;L^2(\R^n))}^2
        \end{split}
        \end{equation}
        for some $C_\rho>0$. So, we may conclude
        \[
        \begin{split}
            &\|u_n(t_0)\|_{L^2(\Omega)}^2+\|(-\Delta)^{s/2}\varphi_n(0)\|_{L^2(\R^n)}^2+2\int_0^{t_0}\|(-\Delta)^{r/2}u_n\|_{L^2(\R^n)}^2\,dt\\
            &\leq C_\eta\left(\|u_{0,n}\|_{L^2(\Omega)}^2+\|(-\Delta)^{r/2}u_{0,n}\|_{L^2(\R^n)}^2+\|u_{1,n}\|_{H^{-s}(\Omega)}^2\right)+C_\rho \|f_n\|_{L^2(0,t_0;H^{-s}(\Omega))}^2\\
            &+C \int_0^{t_0}\left(\|q(t)\|_{L^p(\Omega)}+\rho\right)\left(\|u_n(t)\|_{L^2(\Omega)}^2+\|(-\Delta)^{s/2}\varphi_n(t)\|_{L^2(\R^n)}^2\right)\,dt\\
            &+\eta\left(\|(-\Delta)^{r/2}\varphi_n(0)\|_{L^2(\R^n)}^2+\|(-\Delta)^{s/2}\varphi_n(0)\|_{L^2(\R^n)}^2\right)\\
            &\leq C_\eta\left(\|u_{0,n}\|_{L^2(\Omega)}^2+\|(-\Delta)^{r/2}u_{0,n}\|_{L^2(\R^n)}^2+\|u_{1,n}\|_{H^{-s}(\Omega)}^2\right)+C_\rho \|f_n\|_{L^2(0,t_0;H^{-s}(\Omega))}^2\\
            &+C \int_0^{t_0}\left(\|q(t)\|_{L^p(\Omega)}+\rho\right)\left(\|u_n(t)\|_{L^2(\Omega)}^2+\|(-\Delta)^{s/2}\varphi_n(t)\|_{L^2(\R^n)}^2\right)\,dt\\
            &+\eta (1+C^2_{r,s})\|(-\Delta)^{s/2}\varphi_n(0)\|_{L^2(\R^n)}^2
        \end{split}
        \]
        for any $\eta>0$ and some $C_\eta>0$. In the last estimate, we also used the Poincar\'e inequality between $\widetilde{H}^r(\Omega)$ and $\widetilde{H}^s(\Omega)$ (see \cite[Proposition~2.1]{zimmermann-paperI}), namely
        \begin{equation}
        \label{F:eq: frac poincare ineq}
            \|(-\Delta)^{r/2}\vartheta\|_{L^2(\R^n)}\leq C_{r,s}\|(-\Delta)^{s/2}\vartheta\|_{L^2(\R^n)}
        \end{equation}
        for any $\vartheta\in \widetilde{H}^s(\Omega)$, where $C_{r,s}>0$ only depends on $n,r,s$ and $\Omega$. This implies
        \[
        \begin{split}
            &\|u_n(t_0)\|_{L^2(\Omega)}^2+(1-\eta (1+C_{r,s}^{2}))\|(-\Delta)^{s/2}\varphi_n(0)\|_{L^2(\R^n)}^2+2\int_0^{t_0}\|(-\Delta)^{r/2}u_n\|_{L^2(\R^n)}^2\,dt\\
            &\leq C_\eta\left(\|u_{0,n}\|_{L^2(\Omega)}^2+\|(-\Delta)^{r/2}u_{0,n}\|_{L^2(\R^n)}^2+\|u_{1,n}\|_{H^{-s}(\Omega)}^2\right)+C_\rho \|f_n\|_{L^2(0,t_0;H^{-s}(\Omega))}^2\\
            &+C \int_0^{t_0}\left(\|q(t)\|_{L^p(\Omega)}+\rho\right)\left(\|u_n(t)\|_{L^2(\Omega)}^2+\|(-\Delta)^{s/2}\varphi_n(t)\|_{L^2(\R^n)}^2\right)\,dt.
        \end{split}
        \]
        Next, we again introduce
        \begin{equation}
        \label{F:eq: psi function}
            \psi_n(t)\vcentcolon = \int_0^t u_n(\tau)\,d\tau,
        \end{equation}
        which leads to
        \begin{equation}
         \label{F:eq: almost energy estimate for very weak sol}
        \begin{split}
            &\|u_n(t_0)\|_{L^2(\Omega)}^2+(1-\eta (1+C_{r,s}^{2}))\|(-\Delta)^{s/2}\psi_n(t_0)\|_{L^2(\R^n)}^2+2\int_0^{t_0}\|(-\Delta)^{r/2}u_n\|_{L^2(\R^n)}^2\,dt\\
            &\leq C_\eta\left(\|u_{0,n}\|_{L^2(\Omega)}^2+\|(-\Delta)^{r/2}u_{0,n}\|_{L^2(\R^n)}^2+\|u_{1,n}\|_{H^{-s}(\Omega)}^2\right)+C_\rho \|f_n\|_{L^2(0,t_0;H^{-s}(\Omega))}^2\\
            &+C \int_0^{t_0}\left(\|q(t)\|_{L^p(\Omega)}+\rho\right)\left(\|u_n(t)\|_{L^2(\Omega)}^2+\|(-\Delta)^{s/2}(\psi_n(t_0)-\psi_n(t))\|_{L^2(\R^n)}^2\right)\,dt\\
            &\leq C_\eta\left(\|u_{0,n}\|_{L^2(\Omega)}^2+\|(-\Delta)^{r/2}u_{0,n}\|_{L^2(\R^n)}^2+\|u_{1,n}\|_{H^{-s}(\Omega)}^2\right)+C_\rho \|f_n\|_{L^2(0,t_0;H^{-s}(\Omega))}^2\\
            &+C \int_0^{t_0}\left(\|q(t)\|_{L^p(\Omega)}+\rho\right)\left(\|u_n(t)\|_{L^2(\Omega)}^2+\|(-\Delta)^{s/2}\psi_n(t)\|_{L^2(\R^n)}^2\right)\,dt\\
            &+C_{\ast}(\|q\|_{L^1(0,T;L^p(\Omega))}+\rho T)\|(-\Delta)^{s/2}\psi_n(t_0)\|_{L^2(\R^n)}^2
        \end{split}
        \end{equation}
        for some $C_\ast>0$ only depending on $n,r,s,p$ and $\Omega$. Hence, there holds
        \[
        \begin{split}
             &\|u_n(t_0)\|_{L^2(\Omega)}^2+(1-\eta (1+C_{r,s}^{2})-C_{\ast}(\|q\|_{L^1(0,T;L^p(\Omega))}+\rho T))\|(-\Delta)^{s/2}\psi_n(t_0)\|_{L^2(\R^n)}^2\\
             &+2\int_0^{t_0}\|(-\Delta)^{r/2}u_n\|_{L^2(\R^n)}^2\,dt\\
             &\leq C_\eta\left(\|u_{0,n}\|_{L^2(\Omega)}^2+\|(-\Delta)^{r/2}u_{0,n}\|_{L^2(\R^n)}^2+\|u_{1,n}\|_{H^{-s}(\Omega)}^2\right)+C_\rho \|f_n\|_{L^2(0,t_0;H^{-s}(\Omega))}^2\\
            &+C \int_0^{t_0}\left(\|q(t)\|_{L^p(\Omega)}+\rho\right)\left(\|u_n(t)\|_{L^2(\Omega)}^2+\|(-\Delta)^{s/2}\psi_n(t)\|_{L^2(\R^n)}^2\right)\,dt\\
            &\leq C\left(\|(-\Delta)^{r/2}u_{0}\|_{L^2(\R^n)}^2+\|u_{1}\|_{H^{-s}(\Omega)}^2+ \|f\|_{L^2(0,t_0;H^{-s}(\Omega))}^2\right)\\
            &+C \int_0^{t_0}\left(\|q(t)\|_{L^p(\Omega)}+\rho\right)\left(\|u_n(t)\|_{L^2(\Omega)}^2+\|(-\Delta)^{s/2}\psi_n(t)\|_{L^2(\R^n)}^2\right)\,dt
        \end{split}
        \]
        for sufficiently large $n\in\N$. Now, we choose $0<\eta,\rho\leq 1$ such that $\eta(1+C_{r,s}^{2})\leq 1/4$ and $\rho C_{\ast}T\leq 1/8$. Thus, if condition \eqref{F:eq: smallness condition on q and f} is fulfilled with the constant $C_\ast >0$, then we get 
        \[
            1-\eta (1+C_{r,s}^{2})-C_{\ast}(\|q\|_{L^1(0,T;L^p(\Omega))}+\rho T)\geq 1/2,
        \]
        which ensures that we have the following estimate
        \[
        \begin{split}
            &\|u_n(t)\|_{L^2(\Omega)}^2+\|(-\Delta)^{s/2}\psi_n(t)\|_{L^2(\R^n)}^2+\int_0^{t}\|(-\Delta)^{r/2}u_n\|_{L^2(\R^n)}^2\,d\tau\\
            &\leq C\left(\|(-\Delta)^{r/2}u_{0}\|_{L^2(\R^n)}^2+\|u_{1}\|_{H^{-s}(\Omega)}^2+\|f\|_{L^2(0,t;H^{-s}(\Omega))}^2\right)\\
            &+C\int_0^{t}\left(\|q(\tau)\|_{L^p(\Omega)}+1\right)\left(\|u_n(\tau)\|_{L^2(\Omega)}^2+\|(-\Delta)^{s/2}\psi_n(\tau)\|_{L^2(\R^n)}^2\right)\,d\tau
        \end{split}
        \]
        for all $0\leq t\leq T$ and sufficiently large $n$. Now, using Gronwall's inequality, we deduce
        \begin{equation}
        \label{F:eq: uniform boundedness for very weak sol}
        \begin{split}
            &\|u_n(t)\|_{L^2(\Omega)}^2+\|(-\Delta)^{s/2}\psi_n(t)\|_{L^2(\R^n)}^2+\int_0^{t}\|(-\Delta)^{r/2}u_n\|_{L^2(\R^n)}^2\,d\tau\\
            &\leq C\left(\|(-\Delta)^{r/2}u_{0}\|_{L^2(\R^n)}^2+\|u_{1}\|_{H^{-s}(\Omega)}^2+\|f\|_{L^2(0,t;H^{-s}(\Omega))}^2\right)e^{C\|q\|_{L^1(0,T;L^p(\Omega))}}
        \end{split}
        \end{equation}
        for all $0\leq t\leq T$ and sufficiently large $n$. Therefore, we have demonstrated that $(u_n)_{n\in\N}$ is uniformly bounded in $L^{\infty}(0,T;\widetilde{L}^2(\Omega))$ and $L^2(0,T;\widetilde{H}^r(\Omega))$. Hence, there exists $u\in L^{\infty}(0,T;\widetilde{L}^2(\Omega))\cap L^2(0,T;\widetilde{H}^r(\Omega))$ such that
        \begin{equation}
        \label{F:eq: convergence of un}
            \begin{cases}
                u_n\weakstar u\text{ in }L^{\infty}(0,T;\widetilde{L}^2(\Omega)),\\
                u_n\weak u\text{ in }L^2(0,T;\widetilde{H}^r(\Omega))
            \end{cases}
        \end{equation}
        as $n\to\infty$ (along a suitable subsequence). \\
        
        \emph{Step 4: Passage to the limit and uniqueness.} Recalling that $(u_{0,n},u_{1,n})\to (u_0,u_1)$ in $\widetilde{H}^r(\Omega)\times H^{-s}(\Omega)$ (as $n\to \infty$) and $g\in L^2(0,T;H^{-r}(\Omega))$, we can pass to the limit $n\to\infty$ in \eqref{F:eq: very weak formulation for fn}, which yields
        \begin{equation}
        \label{F:eq: identity for very weak sol}
            \begin{split}
                \int_0^T\langle f,v\rangle\,dt&= \int_0^T \langle u,g\rangle\,dt+\langle u_{0},\partial_t v(0)\rangle_{L^2(\Omega)}-\langle u_{1},v(0)\rangle\\
                &\quad -\langle (-\Delta)^{r/2}u_{0},(-\Delta)^{r/2}v(0)\rangle_{L^2(\R^n)}.
            \end{split}
        \end{equation}
        Furthermore, using \eqref{F:eq: convergence of un} as well as the lower semicontinuity of the norm with respect to weak and weak star convergence, we can pass to the limit in \eqref{F:eq: uniform boundedness for very weak sol} to obtain
        \begin{equation}
        \begin{split}
            &\|u\|_{L^{\infty}(0,t;L^2(\Omega))}^2+\int_0^{t}\|(-\Delta)^{r/2}u\|_{L^2(\R^n)}^2\,dt\\
            &\leq C\left(\|(-\Delta)^{r/2}u_{0}\|_{L^2(\R^n)}^2+\|u_{1}\|_{H^{-s}(\Omega)}^2+\|f\|_{L^2(0,t;H^{-s}(\Omega))}^2\right)e^{C\|q\|_{L^1(0,T;L^p(\Omega))}}
        \end{split}
        \end{equation}
        for any $0\leq t\leq T$.

        Finally, we establish the uniqueness assertion of Theorem \ref{F:thm: ver weak sol}. For this purpose, suppose $\widetilde{u}\in L^{\infty}(0,T;\widetilde{L}^2(\Omega))\cap L^2(0,T;\widetilde{H}^r(\Omega))$ also satisfies \eqref{F:eq: identity for very weak sol}, then $w=u-\widetilde{u}$ fulfills
        \[
        \int_0^{T} \langle w,g\rangle\,dt=0
        \]
        for any $g\in L^2(0,T;H^{-r}(\Omega))$. Since this holds in particular for every $g\in C_c^\infty(\Omega_T)$, we get
$w=0$ in $\mathcal \distr(\Omega_T)$, and hence $w=0$ a.e. in $\Omega_T$. This completes the proof.
    \end{proof}

    It remains to prove the strong convergence assertion used in Step~1.

    \begin{claim}
        \label{F:claim: convergence IV}
            The convergences in \ref{F:convergence 4 very} and \eqref{F:eq: strong convergence of damping term} hold.
        \end{claim}
      
    \begin{proof}[Proof of Claim~\ref{F:claim: convergence IV}.]
            For simplicity we drop in this proof the subscript $n$ for the functions $u_{\eps,n}$, $u_n$ and $(u_{0,n},u_{1,n},f_n)$. 
            
            Let us start by observing that the convergence
            \[
            	u_{\eps}(t)\weak u(t)\text{ in }\widetilde{H}^r(\Omega)\quad\text{and}\ \partial_t u_{\eps}(t)\weakstar \partial_t u(t)\text{ in }H^{-s}(\Omega),
            \]
            for all $0\leq t\leq T$, follow from the trace theorem, \ref{F:convergence 1 very}--\ref{F:convergence 3 very}, and the equation for $u_{\eps}$.

            Next, we show the strong convergence \ref{F:convergence 4 very}. Since the finite energy solution $u$ of \eqref{F:eq: very weak sol approximated} is unique \cite[Theorem~3.22]{zimmermann-paperI}, we know that the whole sequence $(u_\eps)_{\eps>0}$ satisfies \ref{F:convergence 1 very}--\ref{F:convergence 3 very} and not only a subsequence. Furthermore, we may deduce from the estimate \cite[Lemma~3.15]{zimmermann-paperI}  that
            \begin{equation}
            \label{F:eq: weak convergence for strong}
                u_\eps(t)\weak \xi_t\text{ in }\widetilde{H}^s(\Omega)\text{ and }\partial_t u_\eps (t)\weak \eta_t\text{ in }\widetilde{L}^2(\Omega)
            \end{equation}
            for each fixed $0\leq t\leq T$ and some functions $(\xi_t,\eta_t)\in \widetilde{H}^s(\Omega)\times \widetilde{L}^2(\Omega)$. But as observed at the beginning we have
            \[
            u_\eps(t)\weak u(t)\text{ in }\widetilde{H}^r(\Omega)\text{ and }\partial_t u_\eps (t)\weakstar \partial_t u(t)\text{ in }H^{-s}(\Omega)
            \]
            and so by uniqueness of limits we have 
            \begin{equation}
            \label{F:eq: pointwise limits}
                \xi_t=u(t)\text{ and }\eta_t=\partial_t u(t).
            \end{equation}
            Hence, it remains to show that the convergence in \eqref{F:eq: weak convergence for strong} holds strongly. To this end, we introduce the function
            \[
                U_\eps\vcentcolon = u-u_\eps\in C([0,T];\widetilde{H}^s(\Omega))\cap C^1([0,T];\widetilde{L}^2(\Omega))\cap H^1(0,T;\widetilde{H}^r(\Omega))
            \]
            and
            \[
                \psi_\eps(t)\vcentcolon = \|\partial_t U_\eps(t)\|_{L^2(\Omega)}^2+\|(-\Delta)^{s/2}U_\eps(t)\|_{L^2(\R^n)}^2+2 \int_0^t \|(-\Delta)^{r/2}
            \partial_t U_\eps\|_{L^2(\R^n)}^2\,d\tau,
            \]
            for each $0\leq t\leq T$. Using the energy identities \cite[Corollary~3.13]{zimmermann-paperI} for $u_\eps$ and \cite[Proposition~3.20]{zimmermann-paperI} for $u$, we deduce that
            \begin{align}
             \label{F:eq: identity for strong convergence}
             &\\
             &\|\partial_t U_\eps(t)\|_{L^2(\Omega)}^2+\|(-\Delta)^{s/2}U_\eps(t)\|_{L^2(\R^n)}^2+2 \int_0^t \|(-\Delta)^{r/2}
            \partial_t U_\eps\|_{L^2(\R^n)}^2\,d\tau\\
             &+2 \re \int_0^t \langle q U_\eps,\partial_t U_\eps\rangle_{L^2(\Omega)}\,d\tau\\
                 &\leq \|\partial_t U_\eps(t)\|_{L^2(\Omega)}^2+\|(-\Delta)^{s/2}U_\eps(t)\|_{L^2(\R^n)}^2+2 \int_0^t \|(-\Delta)^{r/2}
            \partial_t U_\eps\|_{L^2(\R^n)}^2\,d\tau\\
             &+2\eps\int_0^t\|(-\Delta)^{s/2}\partial_t u_\eps\|_{L^2(\R^n)}^2\,d\tau+2 \re \int_0^t \langle q U_\eps,\partial_t U_\eps\rangle_{L^2(\Omega)}\,d\tau\\
             &=\|\partial_t u(t)\|_{L^2(\Omega)}^2+\|(-\Delta)^{s/2}u(t)\|_{L^2(\R^n)}^2+2 \int_0^t \|(-\Delta)^{r/2}
            \partial_t u\|_{L^2(\R^n)}^2\,d\tau\\
             &+2 \re \int_0^t \langle q u,\partial_t u\rangle_{L^2(\Omega)}\,d\tau\\
             &+\|\partial_t u_\eps(t)\|_{L^2(\Omega)}^2+\|(-\Delta)^{s/2}u_\eps(t)\|_{L^2(\R^n)}^2+2 \int_0^t \|(-\Delta)^{r/2}
            \partial_t u_\eps\|_{L^2(\R^n)}^2\,d\tau\\
             &+2\eps\int_0^t\|(-\Delta)^{s/2}\partial_t u_\eps\|_{L^2(\R^n)}^2\,d\tau+2 \re \int_0^t \langle q u_\eps,\partial_t u_\eps\rangle_{L^2(\Omega)}\,d\tau\\
             &-2\re\langle \partial_t u(t),\partial_t u_\eps(t)\rangle_{L^2(\Omega)}-2\re \langle (-
             \Delta)^{s/2}u(t),(-\Delta)^{s/2}u_\eps(t)\rangle_{L^2(\R^n)}\\
             &-4\re\int_0^t\langle (-\Delta)^{r/2} \partial_t u,(-\Delta)^{r/2} \partial_t u_\eps\rangle_{L^2(\R^n)}\,d\tau\\
             &-2\re\int_0^t\langle qu,\partial_t u_\eps\rangle_{L^2(\Omega)}\,d\tau-2\re\int_0^t\langle qu_\eps,\partial_t u\rangle_{L^2(\Omega)}\,d\tau\\
             &=2(\|u_1\|_{L^2(\Omega)}^2+\|(-\Delta)^{s/2}u_0\|_{L^2(\R^n)}^2)+2\re \int_0^t\langle f,(\partial_t u+\partial_t u_\eps)\rangle\,d\tau\\
             &-2\re\langle \partial_t u(t),\partial_t u_\eps(t)\rangle_{L^2(\Omega)}-2\re \langle (-
             \Delta)^{s/2}u(t),(-\Delta)^{s/2}u_\eps(t)\rangle_{L^2(\R^n)}\\
              &-4\re\int_0^t\langle (-\Delta)^{r/2} \partial_t u,(-\Delta)^{r/2} \partial_t u_\eps\rangle_{L^2(\R^n)}\,d\tau\\
              &-2\re\int_0^t\langle qu,\partial_t u_\eps\rangle_{L^2(\Omega)}\,d\tau-2\re\int_0^t\langle qu_\eps,\partial_t u\rangle_{L^2(\Omega)}\,d\tau\\
            \end{align}
            for all $0\leq t\leq T$. If we introduce the function
            \begin{equation}
                \begin{split}
                    \Psi_\eps(t)\vcentcolon = &2(\|u_1\|_{L^2(\Omega)}^2+\|(-\Delta)^{s/2}u_0\|_{L^2(\R^n)}^2)+2\re \int_0^t\langle f,(\partial_t u+\partial_t u_\eps)\rangle\,d\tau\\
             &-2\re\langle \partial_t u(t),\partial_t u_\eps(t)\rangle_{L^2(\Omega)}-2\re \langle (-
             \Delta)^{s/2}u(t),(-\Delta)^{s/2}u_\eps(t)\rangle_{L^2(\R^n)}\\
             &-4\re\int_0^t\langle (-\Delta)^{r/2} \partial_t u,(-\Delta)^{r/2} \partial_t u_\eps\rangle_{L^2(\R^n)}\,d\tau\\
             &-2\re\int_0^t\langle qu,\partial_t u_\eps\rangle_{L^2(\Omega)}\,d\tau-2\re\int_0^t\langle qu_\eps,\partial_t u\rangle_{L^2(\Omega)}\,d\tau\\
                \end{split} 
            \end{equation}
            and use Lemma~\ref{R:lemma: cond on q} to bound the potential term as
            \begin{equation}
                \left| \int_0^t \langle q U_\eps,\partial_t U_\eps\rangle_{L^2(\Omega)}\,d\tau\right|\leq C\int_0^t \psi_\eps \,d\tau,
            \end{equation}
            then \eqref{F:eq: identity for strong convergence} ensures that
            \begin{equation}
            \label{F:eq: gronwall for strong conv}
                \psi_\eps(t) \leq C\int_0^t \psi_\eps\,d\tau+|\Psi_\eps(t)|.
            \end{equation}
            Next, we apply \eqref{F:eq: weak convergence for strong}, \eqref{F:eq: pointwise limits}, \ref{F:convergence 1 very}--\ref{F:convergence 3 very}, and \cite[Proposition~3.20]{zimmermann-paperI} to get
            \begin{equation}
             \label{F:eq: Psi eps conv}
            \begin{split}
                \Psi_\eps(t)\to &2(\|u_1\|_{L^2(\Omega)}^2+\|(-\Delta)^{s/2}u_0\|_{L^2(\R^n)}^2)+4\re \int_0^t\langle f,\partial_t u\rangle\,d\tau\\
             &-2\| \partial_t u(t)\|^2_{L^2(\Omega)}-2 \| (-
             \Delta)^{s/2}u(t)\|^2_{L^2(\R^n)}\\
             &-4\re\int_0^t\| (-\Delta)^{r/2} \partial_t u\|_{L^2(\R^n)}^2\,dt\\
             &-4\re\int_0^t\langle qu,\partial_t u\rangle_{L^2(\Omega)}\,d\tau\\
             &=0
            \end{split}
            \end{equation}
            as $\eps\to 0$ for any $0\leq t\leq T$. Recalling the energy estimates \cite[Lemma~3.15]{zimmermann-paperI} and  \cite[Proposition~3.16]{zimmermann-paperI}, we see  that $\psi_\eps(t)$ is uniformly bounded in $\eps>0$ and $0\leq t\leq T$. In fact, we have
            \[
            \begin{split}
            	|\psi_{\eps}(t)|&\lesssim \|\partial_t u_\eps(t)\|_{L^2(\Omega)}^2+\|(-\Delta)^{s/2}u_\eps(t)\|_{L^2(\R^n)}^2+2 \int_0^t \|(-\Delta)^{r/2}
            	\partial_t u_\eps\|_{L^2(\R^n)}^2\,d\tau\\
            	&\quad + \|\partial_t u(t)\|_{L^2(\Omega)}^2+\|(-\Delta)^{s/2}u(t)\|_{L^2(\R^n)}^2+2 \int_0^t \|(-\Delta)^{r/2}
            	\partial_t u\|_{L^2(\R^n)}^2\,d\tau\\
            	&\lesssim e^{\alpha T(1+\|q\|_{L^{\infty}(0,T;L^p(\Omega))}^2)}\left(\|u_1\|_{L^2(\Omega)}^2+\|u_0\|_{\widetilde{H}^s(\Omega)}^2+\|f\|_{L^2(0,T;H^{-r}(\Omega))}^2\right)
            \end{split}	
            \]
            for all $\eps>0$ and $0\leq t\leq T$.
            Hence, by \eqref{F:eq: gronwall for strong conv}, Fatou's lemma, and \eqref{F:eq: Psi eps conv}, the function $\psi_0\vcentcolon = \limsup_{\eps\to 0} \psi_\eps$ satisfies 
            \[
                \psi_0\leq C\int_0^t\psi_0\,d\tau.
            \]
            Thus, Gronwall's inequality ensures $\psi_0=0$, which proves
            \begin{equation}
            \label{F:eq: strong convergence}
                \begin{split}
                    u_\eps(t)&\to u(t)\text{ in }\widetilde{H}^s(\Omega),\\
                    \partial_t u_\eps(t)&\to \partial_t u(t)\text{ in }\widetilde{L}^2(\Omega),\\
                    \partial_t u_\eps&\to \partial_t u\text{ in }L^2(0,T;\widetilde{H}^r(\Omega))
                \end{split}
            \end{equation}
            as $\eps\to 0$ for all $0\leq t\leq T$. This establishes property \ref{F:convergence 4 very} and the convergence \eqref{F:eq: strong convergence of damping term}.
       \end{proof}

    \subsection{Very weak theory for equations with \texorpdfstring{$(X^1_{s,r}(\Omega_T))'$}{space-time Sobolev}-sources}
    \label{F:sec: rough sources constant}

    In this subsection we extend the previous construction to initial data $(u_0,u_1)$ and sources $f$ satisfying
    \begin{equation}
    \label{F:eq: rough data}
(u_0,u_1)\in
\begin{cases}
\widetilde H^{2r-s}(\Omega)\times H^{-s}(\Omega),
& s<2r,\\[2mm]
\widetilde L^2(\Omega)\times H^{-s}(\Omega),
& 2r\leq s,
\end{cases}
\qquad f\in (X^1_{s,r}(\Omega_T))'.
\end{equation}
In this case, we use the following definition of very weak solutions.
\begin{definition}[Very weak solutions for energy sources]
\label{F:def: very weak sols for energy sources}
    Let $(u_0,u_1,f)$ satisfy \eqref{F:eq: rough data}.
A function
\[
u\in
L^\infty(0,T;\widetilde L^2(\Omega))
\cap
L^2(0,T;\widetilde H^r(\Omega))
\]
is called a very weak solution of \eqref{F:eq: very weak problem constant coeff} if
for every $g\in L^2(0,T;H^{-r}(\Omega))$, with $v\in  X^1_{s,r,T}(\Omega_T)$ denoting the finite energy
solution of the adjoint terminal-value problem \eqref{F:eq: adjoint problem very weak}, one has the transposition identity
\begin{equation}
\label{F:eq: identity for very weak sol rough}
\begin{split}
&\langle f,v\rangle_{(X^1_{s,r}(\Omega_T))'\times X^1_{s,r}(\Omega_T)}
=
\int_0^T\langle u,g\rangle_{H^{-r}(\Omega)\times \widetilde H^r(\Omega)}\,dt
+
\langle u_0,\partial_t v(0)\rangle_{L^2(\Omega)}\\
&\quad
-
\langle u_1,v(0)\rangle_{H^{-s}(\Omega)\times\widetilde H^s(\Omega)}-
\langle
(-\Delta)^{(2r-s)/2}u_0,
(-\Delta)^{s/2}v(0)
\rangle_{L^2(\R^n)}
\end{split}
\end{equation}
when $s<2r$, and
\begin{equation}
\label{F:eq: identity for very weak sol rough 2}
\begin{split}
&\langle f,v\rangle_{(X^1_{s,r}(\Omega_T))'\times X^1_{s,r}(\Omega_T)}
=
\int_0^T\langle u,g\rangle_{H^{-r}(\Omega)\times \widetilde H^r(\Omega)}\,dt
+
\langle u_0,\partial_t v(0)\rangle_{L^2(\Omega)}
\\
&\quad
-
\langle u_1,v(0)\rangle_{H^{-s}(\Omega)\times\widetilde H^s(\Omega)}-
\langle u_0,(-\Delta)^r v(0)\rangle_{L^2(\Omega)}
\end{split}
\end{equation}
when $2r\leq s$.
\end{definition}

    \begin{theorem}[Very weak solutions with sources in the energy space]
\label{F:thm: very weak sol with rough sources}
There exists a constant
$C_\ast>0$, depending only on $n,r,s,p$, and $\Omega$, such that if
$q\in \mathcal{Q}_s^p$ satisfies
\begin{equation}
\label{F:eq: smallness condition on q rough}
    \|q\|_{L^1(0,T;L^p(\Omega))}
    \leq \frac{1}{8C_\ast},
\end{equation}
then the following assertions hold.

For all data $(u_0,u_1,f)$ satisfying \eqref{F:eq: rough data}, there exists a unique very weak solution
\[
u\in
L^\infty(0,T;\widetilde L^2(\Omega))
\cap
L^2(0,T;\widetilde H^r(\Omega))
\]
of \eqref{F:eq: very weak problem constant coeff} in the sense of Definition~\ref{F:def: very weak sols for energy sources}.

Moreover, for a.e. $0<t<T$, the solution satisfies
\begin{equation}
\label{F:eq: energy estimate high viscosity}
\begin{split}
&\|u(t)\|_{L^2(\Omega)}^2
+
\int_0^t
\|u(\tau)\|_{\widetilde H^r(\Omega)}^2\,d\tau\\
&\quad \lesssim
\left(
\|u_0\|_{\widetilde H^{2r-s}(\Omega)}^2
+
\|u_1\|_{H^{-s}(\Omega)}^2
+
\|f\|_{(X^1_{s,r}(\Omega_T))'}^2
\right)
e^{C\|q\|_{L^1(0,T;L^p(\Omega))}}
\end{split}
\end{equation}
if $s<2r$, and
\begin{equation}
\label{F:eq: energy estimate small viscosity}
\begin{split}
&\|u(t)\|_{L^2(\Omega)}^2
+
\int_0^t
\|u(\tau)\|_{\widetilde H^r(\Omega)}^2\,d\tau\\
&\quad \lesssim
\left(
\|u_0\|_{L^2(\Omega)}^2
+
\|u_1\|_{H^{-s}(\Omega)}^2
+
\|f\|_{(X^1_{s,r}(\Omega_T))'}^2
\right)
e^{C\|q\|_{L^1(0,T;L^p(\Omega))}}
\end{split}
\end{equation}
if $2r\leq s$.
\end{theorem}

  \begin{remark}
    In the regime $0<r<s<2r$, we have
    \[
        0<2r-s<r,
    \]
    and therefore
    \[
        \widetilde H^r(\Omega)\hookrightarrow \widetilde H^{2r-s}(\Omega).
    \]
    Consequently, the assumption on the initial datum $u_0$ in Theorem~\ref{F:thm: very weak sol with rough sources} is weaker than in Theorem~\ref{F:thm: ver weak sol}. Furthermore, $u_0\in \widetilde H^{2r-s}(\Omega)$, $v\in C([0,T];\widetilde{H}^s(\Omega))$, and the continuity properties of the fractional Laplacian ensure that the last term in \eqref{F:eq: identity for very weak sol rough} is well-defined.
\end{remark}

    \begin{proof}
By Lemma~\ref{R:lemma: space of sources}, there exists a sequence
$(f_n)_{n\in\N}\subset L^2(0,T;\widetilde L^2(\Omega))$ such that
\begin{equation}
\label{F:eq: convergence of source rough}
    f_n\to f
    \quad\text{in }(X^1_{s,r}(\Omega_T))'
\end{equation}
as $n\to\infty$. Moreover, by standard approximation, we may choose
\[
    (u_{0,n},u_{1,n})\in \widetilde H^s(\Omega)\times\widetilde L^2(\Omega),
    \qquad n\in\N,
\]
such that
\begin{equation}
\label{F:eq: convergence of initial conds rough}
(u_{0,n},u_{1,n})\to (u_0,u_1)
\quad\text{in}\quad
\begin{cases}
\widetilde H^{2r-s}(\Omega)\times H^{-s}(\Omega),
& s<2r,\\
\widetilde L^2(\Omega)\times H^{-s}(\Omega),
& 2r\leq s.
\end{cases}
\end{equation}

Let $u_n$ be the finite energy solution of
\eqref{F:eq: very weak sol approximated} with data $(u_{0,n},u_{1,n},f_n)$, and let
$v$ be the finite energy solution of the adjoint problem with source
$g\in L^2(0,T;H^{-r}(\Omega))$. By the transposition identity obtained in Step 2 of the proof of
Theorem~\ref{F:thm: ver weak sol}, we have
\begin{equation}
\label{F:eq: very weak formulation for fn rough first}
\begin{split}
\int_0^T\langle f_n,v\rangle\,dt
&=
\int_0^T \langle u_n,g\rangle\,dt
+\langle u_{0,n},\partial_t v(0)\rangle_{L^2(\Omega)}
-\langle u_{1,n},v(0)\rangle
\\
&\quad
-\langle (-\Delta)^{r/2}u_{0,n},
(-\Delta)^{r/2}v(0)\rangle_{L^2(\R^n)}.
\end{split}
\end{equation}

If $s<2r$, then we write this last term as
\[
\langle (-\Delta)^{r/2}u_{0,n},(-\Delta)^{r/2}v(0)\rangle_{L^2(\R^n)}
=
\langle (-\Delta)^{(2r-s)/2}u_{0,n},
(-\Delta)^{s/2}v(0)\rangle_{L^2(\R^n)}.
\]
If $2r\leq s$, then we instead write
\[
\langle (-\Delta)^{r/2}u_{0,n},(-\Delta)^{r/2}v(0)\rangle_{L^2(\R^n)}
=
\langle u_{0,n},(-\Delta)^r v(0)\rangle_{L^2(\R^n)}.
\]

We next derive estimates uniform in $n\in \N$. Let $\varphi_n\in X^1_{s,r,T}(\Omega_T)$ be the function defined in \eqref{F:eq: test function very weak sol} for some fixed $0<t_0<T$.

Arguing as in Step 3 of the proof of Theorem~\ref{F:thm: ver weak sol}, we obtain
\begin{equation}
\label{F:eq: useful identity rough sources}
\begin{split}
&\|u_n(t_0)\|_{L^2(\Omega)}^2
+\|(-\Delta)^{s/2}\varphi_n(0)\|_{L^2(\R^n)}^2
+2\int_0^{t_0}
\|(-\Delta)^{r/2}u_n\|_{L^2(\R^n)}^2\,dt
\\
&=
\|u_{0,n}\|_{L^2(\Omega)}^2
-2\re\langle (-\Delta)^{r/2}u_{0,n},
(-\Delta)^{r/2}\varphi_n(0)\rangle_{L^2(\R^n)}
\\
&\quad
-2\re\langle u_{1,n},\varphi_n(0)\rangle
+2\re\int_0^{t_0}\langle qu_n,\varphi_n\rangle_{L^2(\Omega)}\,dt
-2\re\langle f_n,\varphi_n\rangle_{(X^1_{s,r}(\Omega_T))'\times X^1_{s,r}(\Omega_T)} .
\end{split}
\end{equation}
For every $\rho>0$,
\begin{equation}
\label{F:eq: control source term rough}
\begin{split}
&|\langle f_n,\varphi_n\rangle_{(X^1_{s,r}(\Omega_T))'\times X^1_{s,r}(\Omega_T)} |
\leq
\|f_n\|_{(X^1_{s,r}(\Omega_T))'}
\|\varphi_n\|_{X^1_{s,r}(\Omega_T)}
\\
&\quad\leq
C_\rho\|f_n\|_{(X^1_{s,r}(\Omega_T))'}^2\\
&\quad 
+\rho\left(
\|(-\Delta)^{s/2}\varphi_n\|_{L^2(0,t_0;L^2(\R^n))}^2
+
\|(-\Delta)^{r/2}u_n\|_{L^2(0,t_0;L^2(\R^n))}^2
\right).
\end{split}
\end{equation}

We distinguish the low and high viscosity regimes.

\medskip
\noindent\textbf{Case 1: $s<2r$.}
Using
\[
\langle (-\Delta)^{r/2}u_{0,n},(-\Delta)^{r/2}\varphi_n(0)\rangle
=
\langle (-\Delta)^{(2r-s)/2}u_{0,n},
(-\Delta)^{s/2}\varphi_n(0)\rangle,
\]
we get, for every $\theta>0$,
\[
\begin{split}
&2\left|
\langle (-\Delta)^{r/2}u_{0,n},
(-\Delta)^{r/2}\varphi_n(0)\rangle
\right|\\
&\quad\leq
C_\theta
\|(-\Delta)^{(2r-s)/2}u_{0,n}\|_{L^2(\R^n)}^2
+
\theta
\|(-\Delta)^{s/2}\varphi_n(0)\|_{L^2(\R^n)}^2.
\end{split}
\]
Together with \eqref{F:eq: control source term rough}, the estimate for the potential term (Lemma~\ref{R:lemma: cond on q}), and the same argument as in Step 3 in the proof of Theorem~\ref{F:thm: ver weak sol}, this yields
\[
\begin{split}
&\|u_n(t_0)\|_{L^2(\Omega)}^2
+\|(-\Delta)^{s/2}\psi_n(t_0)\|_{L^2(\R^n)}^2
+\int_0^{t_0}\|(-\Delta)^{r/2}u_n\|_{L^2(\R^n)}^2\,dt
\\
&\leq
C\left(
\|u_0\|_{\widetilde H^{2r-s}(\Omega)}^2
+
\|u_1\|_{H^{-s}(\Omega)}^2
+
\|f\|_{(X^1_{s,r}(\Omega_T))'}^2
\right)
\\
&\quad
+C\int_0^{t_0}
(\|q(t)\|_{L^p(\Omega)}+1)
\left(
\|u_n(t)\|_{L^2(\Omega)}^2
+
\|(-\Delta)^{s/2}\psi_n(t)\|_{L^2(\R^n)}^2
\right)\,dt,
\end{split}
\]
where $\psi_n$ is the function defined in \eqref{F:eq: psi function}.

Here we have chosen $\theta,\rho>0$ sufficiently small and used the smallness
condition \eqref{F:eq: smallness condition on q rough}. Gronwall's inequality gives
\begin{equation}
\label{F:eq: energy inequality rough high viscosity n}
\begin{split}
&\|u_n(t_0)\|_{L^2(\Omega)}^2
+\int_0^{t_0}\|(-\Delta)^{r/2}u_n\|_{L^2(\R^n)}^2\,dt
\\
&\lesssim
\left(
\|u_0\|_{\widetilde H^{2r-s}(\Omega)}^2
+
\|u_1\|_{H^{-s}(\Omega)}^2
+
\|f\|_{(X^1_{s,r}(\Omega_T))'}^2
\right)
e^{C\|q\|_{L^1(0,T;L^p(\Omega))}} .
\end{split}
\end{equation}

\medskip
\noindent\textbf{Case 2: $2r\leq s$.}
In this case,
\[
\langle (-\Delta)^{r/2}u_{0,n},(-\Delta)^{r/2}\varphi_n(0)\rangle
=
\langle u_{0,n},(-\Delta)^r\varphi_n(0)\rangle.
\]
Since $2r\leq s$, the fractional Poincaré inequality gives
\[
\|(-\Delta)^r\varphi_n(0)\|_{L^2(\R^n)}
\lesssim
\|(-\Delta)^{s/2}\varphi_n(0)\|_{L^2(\R^n)}.
\]
Hence, for every $\theta>0$,
\[
2\left|
\langle u_{0,n},(-\Delta)^r\varphi_n(0)\rangle
\right|
\leq
C_\theta\|u_{0,n}\|_{L^2(\Omega)}^2
+
\theta
\|(-\Delta)^{s/2}\varphi_n(0)\|_{L^2(\R^n)}^2.
\]
Proceeding exactly as in Case 1, we obtain
\begin{equation}
\label{F:eq: energy inequality rough low viscosity n}
\begin{split}
&\|u_n(t_0)\|_{L^2(\Omega)}^2
+\int_0^{t_0}\|(-\Delta)^{r/2}u_n\|_{L^2(\R^n)}^2\,dt
\\
&\lesssim
\left(
\|u_0\|_{L^2(\Omega)}^2
+
\|u_1\|_{H^{-s}(\Omega)}^2
+
\|f\|_{(X^1_{s,r}(\Omega_T))'}^2
\right)
e^{C\|q\|_{L^1(0,T;L^p(\Omega))}} .
\end{split}
\end{equation}

In both cases, $(u_n)_{n\in\N}$ is uniformly bounded in
\[
L^\infty(0,T;\widetilde L^2(\Omega))
\cap
L^2(0,T;\widetilde H^r(\Omega)).
\]
Thus, after passing to a subsequence,
\[
u_n\weakstar u
\quad\text{in }L^\infty(0,T;\widetilde L^2(\Omega)),
\qquad
u_n\weak u
\quad\text{in }L^2(0,T;\widetilde H^r(\Omega)).
\]
Taking the limit $n\to \infty$ of \eqref{F:eq: very weak formulation for fn rough first}, using
\eqref{F:eq: convergence of source rough} and \eqref{F:eq: convergence of initial conds rough}, gives the asserted very weak formulation. In the case $s<2r$ this gives
\eqref{F:eq: identity for very weak sol rough}, while in the case $2r\leq s$ it gives
\eqref{F:eq: identity for very weak sol rough 2}. The estimates
\eqref{F:eq: energy estimate high viscosity} and
\eqref{F:eq: energy estimate small viscosity} follow from
\eqref{F:eq: energy inequality rough high viscosity n},
\eqref{F:eq: energy inequality rough low viscosity n}, and lower semicontinuity.

Finally, uniqueness follows exactly as in the proof of
Theorem~\ref{F:thm: ver weak sol}: if two very weak solutions are given, their
difference $w$ satisfies
\[
    \int_0^T\langle w,g\rangle\,dt=0
\]
for all $g\in L^2(0,T;H^{-r}(\Omega))$, and hence $w=0$.
\end{proof}

\section{Weakly viscous equations with constant coefficients}
\label{R:sec: Runge for constant coeff}

In this section we establish several Runge approximation properties for weakly viscous nonlocal wave equations in different function spaces. Throughout this section, $\Omega\subset\R^n$ denotes a bounded Lipschitz domain, $W\subset\Omega_e$ is a fixed nonempty open measurement set, $T>0$ is a finite time-horizon, $0\leq r<s$, and $q\in \mathcal{Q}_s^p$.

We define the associated \emph{Runge set} by
\begin{equation}
\label{R:eq: Runge set}
\mathscr R_W
\vcentcolon=
\{u_\varphi-\varphi\,;\,\varphi\in C_c^\infty(W_T)\},
\end{equation}
where $u_\varphi$ denotes the unique finite energy solution of
\begin{equation}
\label{R:eq: Wave eq for Runge}
\begin{cases}
(\partial_t^2+(-\Delta)^r\partial_t+(-\Delta)^s+q)u=0
&\text{in }\Omega_T,\\
u=\varphi
&\text{in }(\Omega_e)_T,\\
u(0)=0,\ \partial_tu(0)=0
&\text{in }\Omega.
\end{cases}
\end{equation}
We refer to Remark~\ref{R:Remark: non-homogeneous problems} for the well-posedness of nonhomogeneous exterior value problems.

We recall that in \cite{zimmermann-paperI}, we showed that the Runge set is dense in $L^2(0,T;\widetilde{H}^r(\Omega))$. Here, we first prove that, under the additional smallness condition \eqref{R:eq: smallness condition on q rough}, the Runge set is dense in $L^2(0,T;\widetilde H^s(\Omega))$,
and afterwards extend it to the full energy space
\[
X^1_{s,r,0}(\Omega_T)
=
L^2(0,T;\widetilde H^s(\Omega))
\cap
H^1_0(0,T;\widetilde H^r(\Omega)).
\]
The latter result yields a genuine space--time Runge approximation property and constitutes one of the principal results of this article.

We remark that the Runge approximation in $L^2(0,T;\widetilde H^s(\Omega))$
remain valid for arbitrary fixed initial data by translation with the corresponding homogeneous solution; see \cite{Semilinear-nonlocal-wave-eq}. The analogous statement for the energy space $X^1_{s,r,0}(\Omega_T)$ is more delicate due to the imposed trace conditions and leads naturally to an affine approximation problem.

    \subsection{Runge approximation in \texorpdfstring{$L^2(0,T;\widetilde{H}^s(\Omega))$}{space sobolev}}
    \label{R:sec: runge in Hs}

    The proof of the following Runge approximation theorem combines unique continuation properties with a time-regularization argument for very weak solutions.
    
    \begin{theorem}[Runge approximation in $L^2(0,T;\widetilde{H}^s(\Omega))$]
    \label{R:thm: Runge approx in Hs}
        Assume that $q\in\mathcal{Q}_s$ satisfies the smallness condition \eqref{R:eq: smallness condition on q rough}. Let the exponents $r,s$ satisfy one of the following conditions:
\begin{enumerate}[(i)]
\item\label{R:cond 1 exp Runge}
$s\in\N$, $r\in\R_+\setminus\N$;

\item\label{R:cond 2 exp Runge}
$s\in\R_+\setminus\N$, $r\in\N_0$;

\item\label{R:cond 3 exp Runge}
Assume $n\geq2$.
$r,s\in\R_+\setminus\N$ and
\begin{equation}
\begin{cases}
s-r\notin\mathbb N,
&
\quad \text{if } n\in2\mathbb N,
\\[1mm]
s-r\notin\frac12\mathbb N,
& 
\quad \text{if } n\in 2\mathbb N+1.
\end{cases}
\end{equation}
\end{enumerate}
		Then $\mathscr{R}_W$ is dense in $L^2(0,T;\widetilde{H}^s(\Omega))$.
    \end{theorem}

     For the proof of Theorem~\ref{R:thm: Runge approx in Hs} we will need the following lemma for time-independent very weak  solutions.

    \begin{lemma}[Time-independent very weak solutions]
    \label{R:lemma: Time-independent very weak solutions}
        Let $Q\in \mathcal{Q}_s$ and $F\in L^2(0,T;H^{-s}(\Omega))$. Suppose that $V\in \widetilde{H}^r(\Omega)$ is a time-independent very weak solution of 
        \begin{equation}
        \label{R:eq: time indep sol PDE}
            \begin{cases}
				(\partial_t^2 +(-\Delta)^r\partial_t+(-\Delta)^s+Q)V = F&\text{ in }\Omega_T,\\
				V=0  &\text{ in }(\Omega_e)_T,\\
				V(0)=0, \, \partial_t V(0)=0 &\text{ in }\Omega.
			\end{cases}
        \end{equation}
        Then $V = 0$ in $\R^n$.
    \end{lemma}

    \begin{remark}
        In the proof of Theorem~\ref{R:thm: Runge approx in X} it is shown that Lemma~\ref{R:lemma: Time-independent very weak solutions} remains true, when $F\in (X^1_{s,r}(\Omega_T))'$.
    \end{remark}

    \begin{proof}
          Since $V$ is a very weak solution of \eqref{R:eq: time indep sol PDE}, we deduce that there holds
                \begin{equation}
                \label{R:eq: identity very weak sol for V0}
                    \int_0^T \langle F,\omega\rangle\,dt=\int_0^T \langle V,g\rangle \,dt
                \end{equation}
                for any $g\in L^2(0,T;H^{-r}(\Omega))$, where $\omega$ is the unique finite energy solution of 
                \begin{equation}
                \label{R:eq: test solutions Runge Hs}
                    \begin{cases}
				(\partial_t^2 -(-\Delta)^r\partial_t+(-\Delta)^s+Q)\omega = g&\text{ in }\Omega_T,\\
				\omega=0  &\text{ in }(\Omega_e)_T,\\
				\omega(T)=0, \, \partial_t \omega(T)=0 &\text{ in }\Omega.
			\end{cases}
                \end{equation}
                Firstly, we select $\omega(x,t)=\eta(t)\psi(x)$ with $\eta\in C_c^{\infty}((0,T))$ and $\psi \in C_c^{\infty}(\Omega)$. Note that $\omega$ solves \eqref{R:eq: test solutions Runge Hs} in $\Omega_T$ with 
                \begin{equation}
                    g\vcentcolon = \left.\left(\eta''\psi -\eta'(-\Delta)^r \psi+\eta (-\Delta)^s \psi+Q \eta \psi\right)\right|_{\Omega_T}\in L^2(\Omega_T).
                \end{equation}
                Then, \eqref{R:eq: identity very weak sol for V0} becomes
                \[
                     \int_0^T \langle F,\psi\rangle\eta\,dt=\int_0^T \langle V,((-\Delta)^s+Q)\psi\rangle \eta\,dt,
                \]
                because $\eta\in C_c^{\infty}((0,T))$ and $V$ is independent of $t$. Hence, for a.e. $0<t<T$, the function $V\in \widetilde H^r(\Omega)$ satisfies
                \begin{equation}
                \label{R:eq: distributional eq for V}
                \begin{cases}
                    ((-\Delta)^s+Q (t))V=F (t)&\text{ in }\Omega,\\
                    V=0&\text{ in }\Omega_e.
                \end{cases}
                \end{equation}
                Now, we choose the function $\omega$ in \eqref{R:eq: identity very weak sol for V0} to be of the form $\omega(x,t)=\zeta(t)\psi(x)$, $\zeta\in C_b^{\infty}([0,T])$, $\zeta(T)=\partial_t \zeta(T)=0$ and $\psi\in C_c^{\infty}(\Omega)$. The corresponding source $g$ is given by 
                 \begin{equation}
                    g\vcentcolon = \left.\left(\zeta''\psi -\zeta'(-\Delta)^r \psi+\zeta (-\Delta)^s \psi+Q \zeta \psi\right)\right|_{\Omega_T}\in L^2(\Omega_T)
                \end{equation}
                and, in this case, we obtain
                \[
                \begin{split}
                    \int_0^T \langle F,\psi\rangle\zeta\,dt &=\int_0^T (\zeta''\langle V,\psi\rangle-\zeta'\langle V,(-\Delta)^r\psi\rangle)\,dt\\
                    &\quad +\int_0^T \langle V,((-\Delta)^s+Q)\psi\rangle \zeta\,dt\\
                    &=-\zeta'(0)\langle V,\psi\rangle +\zeta(0)\langle V,(-\Delta)^r\psi\rangle\\
                    &\quad +\int_0^T \langle V,((-\Delta)^s+Q)\psi\rangle \zeta\,dt.
                \end{split}
                \]
                Since $V$ is a distributional solution of \eqref{R:eq: distributional eq for V}, the left-hand side and the last integral on the right-hand side cancel. So, we deduce that
                \[
                \zeta'(0)\langle V,\psi\rangle=\zeta(0)\langle V,(-\Delta)^r\psi\rangle.
                \]
                Choosing $\zeta$ so that $\zeta'(0)\neq 0$ and $\zeta(0)=0$, we find $V=0$ in $\Omega$. This finishes the proof.
    \end{proof}

    \begin{proof}[Proof of Theorem~\ref{R:thm: Runge approx in Hs}]
        As in the proof of \cite[Theorem~5.1]{zimmermann-paperI}, the Hahn--Banach theorem reduces the proof to showing that if $f\in L^2(0,T;H^{-s}(\Omega))$ vanishes on $\mathscr{R}_W$, then $f$ vanishes identically. 
        
        For $\varphi\in C_c^\infty(W_T)$, let
        \[
            w=(u_\varphi-\varphi)^\star.
        \]
        Then $w$ solves the adjoint problem
        \begin{equation}
        \label{R:eq: time reversed sol}
			\begin{cases}
				(\partial_t^2 -(-\Delta)^r\partial_t+(-\Delta)^s+q^\star)w = -(-\Delta)^s\varphi^\star+(-\Delta)^r\partial_t \varphi^\star&\text{ in }\Omega_T,\\
				w=0  &\text{ in }(\Omega_e)_T,\\
				w(T)=0, \, \partial_t w(T)=0 &\text{ in }\Omega.
			\end{cases}
		\end{equation}
		Next, let $v$ be the unique very weak solution to 
		\begin{equation}
        \label{R:eq: very weak sol for runge}
			\begin{cases}
				(\partial_t^2 +(-\Delta)^r\partial_t+(-\Delta)^s+q^\star)v = f^\star&\text{ in }\Omega_T,\\
				v=0  &\text{ in }(\Omega_e)_T,\\
				v(0)=0, \, \partial_t v(0)=0 &\text{ in }\Omega;
			\end{cases}
		\end{equation}
        see Theorem~\ref{R:thm: very weak sol}.
		By a change of variables as well as using the facts that $w$ solves \eqref{R:eq: time reversed sol} and $v$ is a very weak solution to \eqref{R:eq: very weak sol for runge}, we have
		\begin{equation}
        \label{R:eq: identity for Runge}
		\begin{split}
			0&=\int_0^T\left\langle f,u_{\varphi}-\varphi\right\rangle\,dt \\
            &= \int_0^T\left\langle f^\star,w\right\rangle\,dt\\
            &=\int_0^T\langle v,(-\Delta)^r\partial_t\varphi^\star-(-\Delta)^s\varphi^\star\rangle_{L^2(\Omega)}\,dt
		\end{split}
		\end{equation}
		for all $\varphi\in C_c^{\infty}(W_T)$.
        
        To apply \cite[Propsition~4.1]{zimmermann-paperI} and \cite[Proposition~4.2]{zimmermann-paperI} slice-wise in time, we require sufficient spatial Sobolev regularity for the functions $v(t)$, $\partial_t v(t)$.
However, such regularity is not available a priori for very weak solutions.

To overcome this difficulty, we regularize $v$ in the time variable. Let $0<\delta_0\ll T/2$. Let us denote by $(\rho_\delta)_{\delta>0}\subset C_c^{\infty}(\R)$ the standard mollifiers in one space dimension. 
        
        Then for all $\varphi\in C_c^{\infty}(W\times (\delta_0,T-\delta_0))$ and all $0<\delta<\delta_0$, we have 
        \[
            \varphi_\delta\vcentcolon =\rho_\delta\ast \varphi\in C_c^{\infty}(W_T).
        \]
        
        Hence, using \eqref{R:eq: identity for Runge} with $\varphi_\delta^\star$ in place of $\varphi$, we obtain
\begin{equation}
\label{R:eq: orthogonality relation RUnge}
\begin{split}
0
&=
\int_0^T
\langle v,
(-\Delta)^r\partial_t\varphi_\delta
-
(-\Delta)^s\varphi_\delta
\rangle_{L^2(\R^n)}
\,dt
\\
&=
-
\int_0^T
\langle
((-\Delta)^r\partial_t+(-\Delta)^s)v_\delta ,\varphi
\rangle
\,dt .
\end{split}
\end{equation}
      If we choose the test function $\varphi$ of the form $\varphi(x,t)=\xi(t)\psi(x)$ with $\xi\in C_c^{\infty}((\delta_0,T-\delta_0))$ and $\psi\in C_c^{\infty}(W)$, then the previous identity ensures
        \begin{equation}
        \label{R:eq: mollified orthogonality relation}
            (-\Delta)^r \partial_t v_\delta +(-\Delta)^s v_\delta =0\text{ in }W\times (\delta_0,T-\delta_0).
        \end{equation}
        where we set 
        \[
            v_\delta\vcentcolon = \rho_\delta\ast v.
        \]
        Taking into account $v=0$ in $(\Omega_e)_T$, $v_\delta(t)=\rho_\delta\ast v(t)\in L^2(\R^n)$ and $\partial_t v_\delta (t)=\rho'_\delta\ast v(t)\in L^2(\R^n)$ for a.e. $t\in\R$, then we can argue as in the proof of \cite[Theorem~5.1]{zimmermann-paperI} to conclude that
        \begin{equation}
        \label{R:eq: cases 1 for runge}
            \partial_t v_\delta =\partial_t(\rho_\delta\ast v)=0\quad\text{in}\quad \R^n\times (\delta_0,T-\delta_0)
        \end{equation}
        in the case \ref{R:cond 1 exp Runge} and  \begin{equation}
        \label{R:eq: cases 2 and 3 for runge}
            v_\delta=\rho_\delta\ast v=0\quad\text{in}\quad \R^n\times (\delta_0,T-\delta_0)
        \end{equation}
        in the case \ref{R:cond 2 exp Runge} or \ref{R:cond 3 exp Runge} for exponents $r,s$ satisfying $s-r\notin\N$ when $n\in 2\N$ (see \cite[Proposition~4.2]{zimmermann-paperI}). 
        
        Next, let us show that in each case \ref{R:cond 1 exp Runge}, \ref{R:cond 2 exp Runge} or \ref{R:cond 3 exp Runge} we have $v\equiv 0$ in $\R^n_T$.
        \begin{enumerate}[{(R}1)]
            \item\label{R:R1} \textit{Case \ref{R:cond 1 exp Runge}:} Notice that \eqref{R:eq: cases 1 for runge} and the fundamental theorem of calculus imply that $\rho_\delta \ast v (t)=\rho_\delta \ast v (t_0)$ for some fixed $t_0\in (\delta_0,T-\delta_0)$ and all $\delta_0<t<T-\delta_0$. If $t,t_0$ are Lebesgue points, passing to the limit $\delta\to0$ yields $v(t)=v(t_0)$ for a.e. $\delta_0<t<T-\delta_0$. By arbitrariness of $\delta_0>0$, we find that
            \begin{equation}
                v(t)=V\text{ for a.e. }0<t<T,
            \end{equation}
            where $V\vcentcolon = v(t_0)\in \widetilde{H}^r(\Omega)$ for some $t_0\in (0,T)$. Now, Lemma~\ref{R:lemma: Time-independent very weak solutions} yields $V\equiv 0$ in $\R^n$.
             \item\label{R:R23} \textit{Case \ref{R:cond 2 exp Runge} or \ref{R:cond 3 exp Runge} with $s-r\notin\N$ when $n\in 2\N$:} As $\rho_\delta\ast v(t)\to v(t)$ in $L^2(\R^n)$ for a.e. $t\in\R$ along a suitable subsequence as $\delta\to 0$, we deduce from \eqref{R:eq: cases 2 and 3 for runge} that $v(t)=0$ in $\R^n\times (\delta_0,T-\delta_0)$. As $\delta_0>0$ can be chosen arbitrarily small, we may conclude that $v=0$ in $\R^n_T$.
        \end{enumerate} 
       Thus, in all cases, $v\equiv0$ in $\R^n_T$. 
Consequently, by the definition of very weak solutions, $f^\star=0$ in the sense of distributions $\distr(\Omega_T)$, and therefore $f\equiv 0$.

By the Hahn--Banach theorem, this proves the density of
$\mathscr R_W$ in
$L^2(0,T;\widetilde H^s(\Omega))$.
    \end{proof}

     \subsection{Runge approximation in \texorpdfstring{$X^1_{s,r,0}(\Omega_T)$}{space--time sobolev}}
    \label{R:sec: runge in X}

    The previous Runge approximation results establish density in Bochner spaces. We now strengthen these results and prove approximation in the full space--time energy space
\[
X^1_{s,r,0}(\Omega_T)
=
L^2(0,T;\widetilde H^s(\Omega))
\cap
H^1_0(0,T;\widetilde H^r(\Omega)).
\]
This yields simultaneous approximation of both the solution and its time derivative and therefore constitutes a genuine space--time Runge approximation property.

In contrast to the $L^2$-setting, the proof in the energy space is substantially more delicate. The main difficulty stems from the structure of the dual space
\[
(X^1_{s,r,0}(\Omega_T))',
\]
which contains contributions acting on time derivatives and traces. As a consequence, the orthogonality relations arising in the Hahn--Banach argument no longer reduce to a purely bulk equation. 

    \begin{theorem}[{Runge approximation in $X^1_{s,r,0}(\Omega_T)$}]
    \label{R:thm: Runge approx in X}
        Suppose that $q\in\mathcal{Q}_s$ satisfies the smallness condition \eqref{R:eq: smallness condition on q rough}. Assume that the exponents $r,s$ satisfy \ref{R:cond 1 exp Runge}, \ref{R:cond 2 exp Runge}, or \ref{R:cond 3 exp Runge} of Theorem~\ref{R:thm: Runge approx in Hs}. 
        
        Then $\mathscr{R}_W$ is dense in $X^1_{s,r,0}(\Omega_T)$.
    \end{theorem}

\begin{proof}
    For later use, we introduce the time-reversed Runge set
\[
\mathscr R_W^\star
\vcentcolon=
\left\{
(u_\varphi-\varphi)^\star
\,;\,
\varphi\in C_c^\infty(W_T)
\right\}.
\]
Since
\[
(X^1_{s,r,0}(\Omega_T))^\star
=
X^1_{s,r,T}(\Omega_T),
\]
it suffices to prove that $\mathscr R_W^\star$ is dense in
$X^1_{s,r,T}(\Omega_T)$.

We proceed by a Hahn--Banach argument. Let $f\in (X^1_{s,r,T}(\Omega_T))'$ and assume that
\[
\langle f,\psi\rangle =0
\qquad\text{for all }\psi\in \mathscr R_W^\star.
\]
We show that $f=0$.

By Lemma~\ref{R:lemma: space of sources}\,\ref{R:Dual of X}, we may write
\[
f=f_1+f_2
\]
with
\[
f_1\in L^2(0,T;H^{-s}(\Omega)),
\qquad
f_2\in (H^1_T(0,T;\widetilde H^r(\Omega)))'.
\]
Since $H^1_T(0,T;\widetilde H^r(\Omega))$ is a (closed) subspace of
$H^1(0,T;\widetilde H^r(\Omega))$, the Hahn--Banach theorem yields an extension
\[
\widetilde f_2\in (H^1(0,T;\widetilde H^r(\Omega)))'
\]
of $f_2$. Hence
\begin{equation}
\label{R:eq: extension of source final}
\widetilde f \vcentcolon= f_1+\widetilde f_2 \in (X^1_{s,r}(\Omega_T))',
\qquad
\widetilde f|_{X^1_{s,r,T}(\Omega_T)}=f.
\end{equation}

Let $v$ be the unique very weak solution to
\[
\begin{cases}
(\partial_t^2+(-\Delta)^r\partial_t+(-\Delta)^s+q^\star)v=\widetilde f
&\text{in }\Omega_T,\\
v=0&\text{in }(\Omega_e)_T,\\
v(0)=0,\ \partial_t v(0)=0&\text{in }\Omega.
\end{cases}
\]

Let $\varphi\in C_c^\infty(W_T)$. Using $f|_{\mathscr R_W^\star}=0$, $(u_\varphi-\varphi)^\star \in X^1_{s,r,T}(\Omega_T)$,
the extension property
\eqref{R:eq: extension of source final},
the definition of very weak solutions,
and the time-reversal identity
\eqref{R:eq: time reversed sol},
we obtain
\[
\begin{split}
0
&=\langle f,(u_\varphi-\varphi)^\star\rangle\\
&=\langle \widetilde f,(u_\varphi-\varphi)^\star\rangle\\
&=\int_0^T
\big\langle
v,\,-(-\Delta)^s\varphi^\star+(-\Delta)^r\partial_t\varphi^\star
\big\rangle_{L^2(\Omega)}
\,dt.
\end{split}
\]
Thus
\begin{equation}
\label{R:eq: orthogonality relation final}
\int_0^T
\big\langle
v,\,-(-\Delta)^s\varphi^\star+(-\Delta)^r\partial_t\varphi^\star
\big\rangle_{L^2(\Omega)}
\,dt
=0
\qquad\text{for all }\varphi\in C_c^\infty(W_T).
\end{equation}

We now show that $v=0$ in $\R^n_T$ in each case.

\medskip

\begin{itemize}
    \item \textit{Case \ref{R:cond 1 exp Runge}:} As in the proof of Theorem~\ref{R:thm: Runge approx in Hs}, the exterior orthogonality identity \eqref{R:eq: orthogonality relation final} with exterior conditions of the form $\rho_\delta \ast \varphi$, where $\rho_\delta$ is a standard mollifier, imply that
\[
    \partial_t v=0
    \quad\text{in }\mathcal D'(0,T;\widetilde H^r(\Omega)).
\]
Consequently, there exists $V\in\widetilde H^r(\Omega)$ such that
\[
    v(x,t)=V(x)
    \quad\text{for a.e. }(x,t)\in\R^n_T.
\]
We show that $V=0$ using the continuity of $\widetilde f$ on $X^1_{s,r}(\Omega_T)$. This result generalizes Lemma~\ref{R:lemma: Time-independent very weak solutions}.

Fix $\psi\in C_c^\infty(\Omega)$ and choose
$\chi\in C_c^\infty(\R)$ with, $\chi\geq 0$, $\chi=1$ in a neighborhood of $0$.
For sufficiently small $\varepsilon>0$, set
\[
    \eta_\varepsilon(t)=t\chi(t/\varepsilon),
    \qquad
    w_\varepsilon(x,t)=\eta_\varepsilon(t)\psi(x),
    \qquad 0\leq t\leq T.
\]
Then $w_\varepsilon$ vanishes near $t=T$, and
\[
    \eta_\varepsilon(0)=0,
    \qquad
    \eta_\varepsilon'(0)=1.
\]
Moreover,
\[
    \|\eta_\varepsilon\|_{L^2(0,T)}
    \lesssim\varepsilon^{3/2},
    \qquad
    \|\eta_\varepsilon'\|_{L^2(0,T)}
    \lesssim\varepsilon^{1/2}.
\]
Thus
\[
    w_\varepsilon\to0
    \quad\text{in }X^1_{s,r}(\Omega_T),
\]
and hence
\begin{equation}
\label{R:eq: stationary mode source limit}
    \langle\widetilde f,w_\varepsilon\rangle\to0.
\end{equation}

On the other hand, $w_\varepsilon$ is the finite energy solution
of the adjoint problem with zero terminal data and source
\[
    g_\varepsilon
    =
    \eta_\varepsilon''\psi
    -\eta_\varepsilon'(-\Delta)^r\psi
    +\eta_\varepsilon(-\Delta)^s\psi
    +q^\star\eta_\varepsilon\psi
    \in L^2(0,T;H^{-r}(\Omega)).
\]
It is therefore an admissible test function in the transposition
identity. Since $v=V$ in $H^r(\R^n)$, this gives
\[
\begin{split}
    \langle\widetilde f,w_\varepsilon\rangle
    &=
    \int_0^T\langle V,g_\varepsilon\rangle_{L^2(\Omega)}\,dt\\
    &=
    -\langle V,\psi\rangle_{L^2(\Omega)}
    +\int_0^T\eta_\varepsilon(t)
      \langle V,((-\Delta)^s+q^\star(t))\psi
      \rangle_{L^2(\Omega)}\,dt,
\end{split}
\]
where we used
\[
    \int_0^T\eta_\varepsilon''(t)\,dt=-\eta'_\eps(0)=-1,
    \qquad
    \int_0^T\eta_\varepsilon'(t)\,dt=-\eta_\eps (0)=0.
\]
Because $V\in L^2(\Omega)$, $\psi\in C_c^\infty(\Omega)$,
and $q^\star\in L^\infty(0,T;L^p(\Omega))$ with $p\geq2$,
we have
\[
    \operatorname*{ess\,sup}_{0<t<T}
    \left|
    \langle V,((-\Delta)^s+q^\star(t))\psi
    \rangle_{L^2(\Omega)}
    \right|<\infty.
\]
Together with
$\|\eta_\varepsilon\|_{L^1(0,T)}\lesssim\varepsilon^2$,
this implies
\[
    \langle\widetilde f,w_\varepsilon\rangle
    \to
    -\langle V,\psi\rangle_{L^2(\Omega)}.
\]
Comparing with \eqref{R:eq: stationary mode source limit}, we obtain
$\langle V,\psi\rangle_{L^2(\Omega)}=0$.
Since $\psi\in C_c^\infty(\Omega)$ was arbitrary and $V$ vanishes
outside $\Omega$, we conclude that $V=0$ in $\R^n$.
\item \textit{Case \ref{R:cond 2 exp Runge} and \ref{R:cond 3 exp Runge}:} If the exponents satisfy
\ref{R:cond 2 exp Runge} or \ref{R:cond 3 exp Runge}, the argument in the proof of Theorem~\ref{R:thm: Runge approx in Hs} directly yields $v=0$.
\end{itemize}

Finally, let $w_k=\eta\varphi_k$, where $\eta\in C_c^\infty([0,T))$ and $\varphi_k\in \widetilde H^s(\Omega)$ is an eigenfunction of the fractional Laplacian $ (-\Delta)^s $ in $\Omega$ with eigenvalue $\lambda_k$. Then 
\[
    w_k\in X^1_{s,r,T}(\Omega_T)\cap C^{\infty}_c([0,T);\widetilde{H}^s(\Omega)),
\]
and moreover
\[
(\partial_t^2-(-\Delta)^r\partial_t+(-\Delta)^s+q^\star)w_k = g_k
\]
in $\Omega_T$, where
\[
g_k
\vcentcolon =
\left.(\eta''\varphi_k-\eta'(-\Delta)^r\varphi_k+\lambda_k\eta\varphi_k+q^\star \eta\varphi_k)\right|_{\Omega_T}
.
\]
Since $\varphi_k\in \widetilde H^s(\Omega)$, $ (-\Delta)^r\varphi_k\in H^{s-2r}(\R^n)$, and by the assumption on the exponents we have $H^{s-2r}(\R^n)\hookrightarrow H^{-r}(\R^n)$, it follows that
\[
g_k\in L^2(0,T;H^{-r}(\Omega)).
\]
Hence $w_k$ is an admissible test function in the definition of very weak solutions. Since $v=0$, we obtain
\[
\langle \widetilde f,w_k\rangle =0.
\]
By \eqref{R:eq: extension of source final}, this yields $\langle f,w_k\rangle=0 $ for all such $w_k$.

Finally, Lemma~\ref{R:lemma: space of sources}\,\ref{R:density in X T}
shows that the linear span of the functions
\[
w_k=\eta\varphi_k,
\qquad
\eta\in C_c^\infty([0,T)),
\]
is dense in $X^1_{s,r,T}(\Omega_T)$. Since $f\in (X^1_{s,r,T}(\Omega_T))'$ is continuous and
\[
\langle f,w_k\rangle=0
\]
for all such functions, we conclude that $f=0$.
\end{proof}

\begin{remark}[Elimination of stationary modes]
In the case $s\in\N$ and $r\notin\N_0$, the exterior
orthogonality identity initially yields only that the associated very weak solution is independent of time. Nevertheless, this stationary mode must vanish. Indeed, the transposition identity detects its pairing with the initial velocity of an adjoint test function. The functions $w_\varepsilon$ constructed above converge to zero in $X^1_{s,r,T}(\Omega_T)$ while satisfying $\partial_t w_\varepsilon(0)=\psi$. Continuity of the source functional therefore forces $\langle V,\psi\rangle=0$ for every $\psi\in C_c^\infty(\Omega)$.
\end{remark}

  \begin{remark}
        As noticed earlier, for the $L^2$-Runge approximation results, the assumption of vanishing initial data can be removed by translation with the corresponding homogeneous solution. In the space--time energy setting, however, the approximation problem becomes affine due to the imposed trace conditions.
    \end{remark}

\section{An inverse problem for nonlinear velocity damping}
\label{sec:nonlinear-velocity-damping}

In this section, we consider the problem
\begin{equation}
\label{eq:nld-problem}
\begin{cases}
   (\partial_t^2+(-\Delta)^r \partial_t +(-\Delta)^s)u+qN_p(\partial_t u)=0\ &\text{in }\Omega_T,\\
    u=\varphi\ &\text{in }(\Omega_e)_T,\\
    u(0)=\partial_t u(0)=0\ &\text{in }\Omega,
\end{cases}
\end{equation}
where $0<r<s$, $\varphi\in C_c^{\infty}((\Omega_e)_T)$, $q\in L^\infty(\Omega_T)$ is nonnegative, and the nonlinearity $N_p$ is given by
\begin{equation}
\label{eq:nld-range}
N_p(\tau)=|\tau|^{p-2} \tau\quad\text{with}\quad
 2<p\leq 2+4r/n.
\end{equation}
Here, the endpoint is included, so that we consider an energy-critical-or-subcritical range. We remark that throughout this section, all data, solutions, and test functions
are real-valued.

The goal of this section is to study the unique recovery of $q$ from the associated Dirichlet to Neumann map $\Lambda_q$.

\subsection{Nonlinear forward theory}
\label{sec:nonlinear-forward}

We observe that for exterior data $\varphi\in C_c^{\infty}((\Omega_e)_T)$ there is a unique solution of \eqref{eq:nld-problem} with
\[
    u-\varphi \in C([0,T];\widetilde H^{s}(\Omega))\cap C^1([0,T];L^2(\Omega))
       \cap H^1(0,T;\widetilde H^{r}(\Omega)).
\]
Indeed, as in the linear forward theory \cite{zimmermann-paperI}, we obtain this solution by parabolic regularization. We need to use that regularization procedure, because the energy bounds in the weakly viscous regime only control the velocity in $L^2(0,T;\widetilde{H}^r(\Omega))$, whereas the elastic energy controls $U$ in  $L^\infty (0,T;\widetilde{H}^{s}(\Omega))$. For the zero-exterior lift $U=u-\varphi$, we therefore consider the regularized problems
\[
 \partial_t^2U_\varepsilon+\varepsilon (-\Delta)^s\partial_tU_\varepsilon
 +(-\Delta)^r\partial_tU_\varepsilon+(-\Delta)^s U_\varepsilon
 +q N_p(\partial_tU_\varepsilon)=F_\varphi,
 \qquad \varepsilon>0,
\]
with zero exterior and initial data, where
$F_\varphi=-((-\Delta)^r\partial_t\varphi+(-\Delta)^s\varphi)|_\Omega$.
The additional viscous term gives
$\partial_tU_\varepsilon\in L^2(0,T;\widetilde{H}^s(\Omega))$; see \cite[Section~3.1]{zimmermann-paperI}.
The regularized problems are treated by the viscous evolution framework \cite[Definition~3.1]{zimmermann-paperI}. 

For clarity, we record the estimates underlying local solvability of
these regularized problems. Fix $\varepsilon>0$, let $I=(t_0,t_0+\tau)$,
and put $Q_I=\Omega\times I$ and $P=2+4s/n$. We work in the Banach space
\[
 \mathcal E_\varepsilon(I)
 =
 C(\overline I;\widetilde H^s(\Omega))
 \cap C^1(\overline I;L^2(\Omega))
 \cap H^1(I;\widetilde H^s(\Omega)),
\]
equipped with the following energy norm:
\[
\begin{split}
 \|V\|_{\mathcal{E}_\varepsilon (I)}={}&
 \|V\|_{L^\infty(I;\widetilde{H}^s(\Omega))}
 +\|\partial_t V\|_{L^\infty(I;L^2(\Omega))}\\
 &+\|\partial_t V\|_{L^2(I;\widetilde{H}^r(\Omega))}
 +\sqrt{\varepsilon}\|\partial_t V\|_{L^2(I;\widetilde{H}^s(\Omega))}.
\end{split}
\]
The space--time Sobolev inequality, with $s$ in place of $r$ in
Lemma~\ref{lem:nld-embedding}, gives
\begin{equation}
\label{eq:nld-regularized-embedding}
 \|\partial_t V\|_{L^P(Q_I)}\leq C_\varepsilon
                      \|V\|_{\mathcal{E}_\varepsilon(I)}.
\end{equation}
Here and below the constants can depend on $\varepsilon$, but are
uniform for subintervals of $[0,T]$.
If 
\[
    (\partial_t^2+\eps (-\Delta)^s\partial_t+(-\Delta)^r \partial_t +(-\Delta)^s) Z=G\quad\text{in }Q_I
\]
with zero exterior and initial data at $t_0$, velocity testing gives
\[
 \|Z\|_{\mathcal{E}_\varepsilon(I)}^2
 \leq C_\varepsilon\sup_{t\in\overline I}\left|\int_{t_0}^t\int_\Omega G\partial_t Z\,dx\,d\rho\right|
 \leq C_\varepsilon\|G\|_{L^{P'}(Q_I)}
                         \|Z\|_{\mathcal{E}_\varepsilon(I)}.
\]
Thus the linear solution operator satisfies
\begin{equation}\label{eq:nld-regularized-source}
 \|Z\|_{\mathcal{E}_\varepsilon(I)}
 \leq C_\varepsilon\|G\|_{L^{P'}(Q_I)},\qquad P'=P/(P-1).
\end{equation}
For general $G\in L^{P'}(Q_I)$, approximate $G$ in this space
by smooth sources. Applying the preceding estimate to differences
shows that the corresponding solutions form a Cauchy sequence in
$\mathcal E_\varepsilon(I)$.
Their limit defines the linear solution operator on $L^{P'}(Q_I)$
and satisfies the same estimate.
For nonzero initial data, the right-hand side additionally contains
$C_\varepsilon(\|Z(t_0)\|_{\widetilde{H}^s(\Omega)}+\|\partial_t Z(t_0)\|_{L^2(\Omega)})$; cf. \cite[Corollary~3.13]{zimmermann-paperI}.

Since $p\leq2+4r/n<P$, H\"older's inequality and the pointwise
Nemytskii estimate yield
\begin{equation}\label{eq:nld-regularized-lipschitz}
\begin{split}
 &\|q(N_p(\partial_t V)-N_p(\partial_t W))\|_{L^{P'}(Q_I)}\\
 &\quad\leq C_p\|q\|_\infty |Q_I|^{1-p/P}
 (\|\partial_t V\|_{L^P(Q_I)}+\|\partial_t W\|_{L^P(Q_I)})^{p-2}
 \|\partial_t V-\partial_t W\|_{L^P(Q_I)}.
\end{split}
\end{equation}
Define $\mathcal S(V)$ as the solution of the regularized linear
problem with source $F_\varphi-qN_p(\partial_t V)$ and the prescribed initial
data. Combining the preceding estimates, on an energy ball of radius
$R$ we obtain
\[
 \|\mathcal S(V)-\mathcal S(W)\|_{\mathcal{E}_\varepsilon(I)}
 \leq C_\varepsilon\|q\|_\infty
 \tau^{1-p/P}R^{p-2}\|V-W\|_{\mathcal{E}_\varepsilon(I)}.
\]
Choose $R$ to be twice a bound for the linear solution with source
$F_\varphi$ and the given initial data; cf.~\cite{Semilinear-nonlocal-wave-eq}. Taking $\tau$ sufficiently small
makes $\mathcal S$ preserve this ball and be a strict contraction.
This proves local solvability for each fixed $\varepsilon>0$.

Velocity testing for the nonlinear solution yields the energy identity
below with the additional nonnegative term
$\varepsilon\int_0^t \|(-\Delta)^{s/2}\partial_tU_\varepsilon\|_{L^2(\R^n)}^2\,d\rho$.
The resulting bound on the displacement and velocity traces permits
iteration of the local construction to reach $T$; no smallness of the
exterior datum is needed. The local time step may depend on
$\varepsilon$, whereas these dissipative energy bounds are uniform
as $\varepsilon\downarrow0$.

The uniform energy estimates and the regularized equation bound the second order time derivative
$\partial_t^2U_\varepsilon$ in
$L^{p'}(0,T;H^{-m}(\Omega))$ for
$m>\max\{s,n/2\}$.
Therefore, by the compact embedding $\widetilde H^r(\Omega)\hookrightarrow L^2(\Omega)$, the
Aubin--Lions compactness lemma gives, along a subsequence,
strong convergence of $\partial_tU_\varepsilon$ in
$L^2(\Omega_T)$ and convergence almost everywhere. The uniform $L^{p'}(\Omega_T)$ bound on
$qN_p(\partial_tU_\varepsilon)$ identifies its weak limit as $qN_p(U_t)$. Moreover, the weighted regularization estimate (cf.~\cite[Lemma~3.15]{zimmermann-paperI}) gives
\[
 \|\varepsilon (-\Delta)^s\partial_tU_\varepsilon\|_{L^2(0,T;H^{-s}(\Omega))}
 \leq C\sqrt{\varepsilon}\to 0.
\]
The energy argument, with parabolic regularization as in \cite[Section~3.3]{zimmermann-paperI},
gives the energy identity and strong continuity of $U$ in $\widetilde{H}^s(\Omega)$
and of $U_t$ in $L^2(\Omega)$. Finally, the uniqueness argument of \cite[Section~3.4]{zimmermann-paperI} also applies to our setting, when one takes into account the monotonicity of $N_p$. In particular, on writing $U=u-\varphi$ and
\[
    F_\varphi=-((-\Delta)^r\partial_t\varphi+(-\Delta)^s\varphi)|_\Omega \in L^2(0,T;L^2(\Omega)),
\]
we have the energy identity
\begin{equation}\label{eq:nld-forward-energy}
\begin{split}
     &E_U(t)+\int_0^t \|(-\Delta)^{r/2}\partial_t U\|_{L^2(\R^n)}^2\,d\rho
       +\int_{\Omega\times(0,t)}q|\partial_t U|^p\,dxdt\\
      & \quad
 =\int_0^t(F_\varphi,\partial_t U)_{L^2(\Omega)}\,d\rho,
\end{split}
\end{equation}
where
\[
 E_U(t)=\tfrac12\bigl(\|\partial_t U(t)\|_{L^2(\Omega)}^2+\|(-\Delta)^{s/2}U(t)\|_{L^2(\R^n)}^2\bigr).
\]
Since the nonnegative nonlinear term can be discarded in the energy
estimate, we see that $u_\delta=u_{q,\delta\varphi}$ obeys
\begin{equation}\label{eq:nld-scaled-bound}
 \|u_\delta-\delta\varphi\|_{\mathcal{E}}\leq C_\varphi\delta
 \qquad(0<\delta\leq1),
\end{equation}
where
\[
    \|\cdot\|_{\mathcal{E}}=\|\cdot\|_{L^{\infty}(0,T;\widetilde{H}^s(\Omega))}+\|\cdot\|_{W^{1,\infty}(0,T;L^2(\Omega))}+\|\cdot\|_{H^1(0,T;\widetilde{H}^r(\Omega))}.
\]
The same testing argument for differences uses
$(N_p(a)-N_p(b))(a-b)\geq0$.
The same approximation and time-regularization argument justifies
velocity testing for the differences used below.

\subsection{Sobolev embeddings and source spaces}

We record here the following space--time Sobolev embedding.
\begin{lemma}[Space--time Sobolev embeddings]\label{lem:nld-embedding}
If $2\leq m\leq2+4r/n$, then
\[
 X^1_{s,r,\tau}(\Omega_T)\hookrightarrow L^m(\Omega_T),\qquad \tau=0,T.
\]
Furthermore, with $p_*=2+4r/n$,
\begin{equation}\label{eq:nld-velocity-embedding}
 L^\infty(0,T;L^2(\Omega))\cap L^2(0,T;\widetilde H^{r}(\Omega))
 \hookrightarrow L^{p_*}(\Omega_T).
\end{equation}
\end{lemma}
\begin{proof}
The Hilbert-valued Sobolev embedding gives
\[
H^1(0,T;\widetilde H^{r}(\Omega))\hookrightarrow C([0,T];\widetilde H^{r}(\Omega)).
\]
If $r<n/2$, spatial Sobolev embedding applies up to
$2n/(n-2r)$, and
\[
 \frac{2n}{n-2r}-\left(2+\frac{4r}{n}\right)
 =\frac{8r^2}{n(n-2r)}>0.
\]
For $r=n/2$ every finite spatial exponent is allowed \cite{Ozawa}; for $r>n/2$
one has the spatial $L^\infty$ embedding. Thus
$\|w\|_{L^m(\Omega_T)}\leq C T^{1/m}\|w\|_{H^1(0,T;\widetilde H^{r}(\Omega))}$.
This embedding in fact holds throughout the corresponding spatial
Sobolev range, without the restriction $m\leq3$.

For the second assertion, the Gagliardo--Nirenberg inequality on the
zero extension (or interpolation theory \cite{Interpolation-spaces-Bergh-Lofstrom}) gives
\[
 \|b(t)\|_{L^{p_*}(\Omega)}
 \leq C\|b(t)\|_{H^r(\R^n)}^\theta
          \|b(t)\|_{L^2(\Omega)}^{1-\theta},
 \qquad \theta=\frac{n}{n+2r},\quad \theta p_*=2.
\]
Integrating its $p_*$th power proves \eqref{eq:nld-velocity-embedding}.
Lower exponents follow since $|\Omega_T|<\infty$.
\end{proof}

Consequently, the velocities of the forward energy solutions belong
to $L^p(\Omega_T)$ in the indicated range. This assertion applies to their
zero-exterior lifts; the exterior data and all their time derivatives
vanish in $\Omega_T$. In particular,
\[
 qN_p(\partial_t u)\in L^{p'}(\Omega_T),\qquad p'=p/(p-1),
 \qquad L^{p'}(\Omega_T)\hookrightarrow (X^1_{s,r,T}(\Omega_T))'.
\]

\subsection{The linear source-to-boundary identity and its sign}
Now, consider the solution $v_\varphi$ of
\begin{equation}
\label{eq:v-phi}
(\partial_t^2+(-\Delta)^r \partial_t +(-\Delta)^s)v_\varphi=0\quad\text{in }\Omega_T
\end{equation}
with exterior value $\varphi\in C_c^{\infty}((\Omega_e)_T)$
and zero initial data. Let $w_\zeta$, $\zeta\in C_c^{\infty}((\Omega_e)_T)$ solve the adjoint problem
\begin{equation}
\label{eq:adjoint-eq}
\begin{cases}
     (\partial_t^2-(-\Delta)^r \partial_t +(-\Delta)^s) w_\zeta=0,\quad &\text{in }\Omega_T\\
 w_\zeta=\zeta,\quad &\text{in }(\Omega_e)_T,\\
 w_\zeta(T)=(w_\zeta)_t(T)=0,\quad&\text{in } \Omega.
\end{cases}
\end{equation}
Time reversal of the linear forward problem constructs this adjoint
solution \cite{zimmermann-paperI}. The lifts $v_\varphi-\varphi$ and $w_\zeta-\zeta$ belong,
respectively, to $X^1_{s,r,0}(\Omega_T)$ and $X^1_{s,r,T}(\Omega_T)$.
The restrictions $w_\zeta|_{\Omega_T}$ therefore belong to $L^p(\Omega_T)$.

Throughout this section, the Dirichlet to Neumann map associated with the nonlinear problem \eqref{eq:nld-problem} is
\begin{equation}\label{eq:nld-DN}
 \langle\Lambda_q(\varphi),\zeta\rangle
 =\int_0^T\langle (-\Delta)^{r/2}\partial_t u,(-\Delta)^{r/2}\zeta\rangle_{L^2(\R^n)}+\langle (-\Delta)^{s/2}u,(-\Delta)^{s/2}\zeta\rangle_{L^2(\R^n)}\,dt,
\end{equation}
for all $\zeta\in C_c^\infty((\Omega_e)_T)$, where $u$ is the unique solution of \eqref{eq:nld-problem}.
Thus $\Lambda_q(\varphi)=((-\Delta)^r\partial_t u+(-\Delta)^s u)|_{(\Omega_e)_T}$ as a
distribution in $(\Omega_e)_T$.

\begin{lemma}[Linear source-to-boundary Green identity]\label{lem:nld-green}
Let 
\[
    H\in C([0,T];\widetilde H^{s}(\Omega))\cap C^1([0,T];L^2(\Omega))
       \cap H^1(0,T;\widetilde H^{r}(\Omega))
\]
have zero exterior and initial data, and satisfy
\[
    (\partial_t^2+(-\Delta)^r\partial_t +(-\Delta)^s) H=f\quad\text{in }\Omega_T
\]
with $f\in L^{p'}(\Omega_T)$. Then
\begin{equation}\label{eq:nld-green}
\begin{split}
   & \int_0^T\langle (-\Delta)^{r/2}\partial_t H,(-\Delta)^{r/2}\zeta\rangle_{L^2(\R^n)}\,dt\\
   &\quad +\int_0^T \langle (-\Delta)^{s/2}H,(-\Delta)^{s/2}\zeta\rangle_{L^2(\R^n)}\,dt
 =-\int_{\Omega_T}f w_\zeta\,dx\,dt.
\end{split}
\end{equation}
\end{lemma}
\begin{proof}
    The result follows by approximation from the weak formulation, Lemma~\ref{lem:nld-embedding}, and the adjoint equation \eqref{eq:adjoint-eq}.
\end{proof}

\subsection{Small-amplitude response proposition}
\begin{proposition}
\label{prop:nld-response}
Let $q\in L^\infty(\Omega_T)$, $q\geq0$, and $2<p\leq2+4r/n$. For any fixed smooth exterior data $\varphi,\zeta\in C_c^{\infty}((\Omega_e)_T)$, one has
\begin{equation}\label{eq:nld-response}
 \lim_{\delta\downarrow0}\delta^{-(p-1)}
 \langle\Lambda_q(\delta\varphi)-\delta\Lambda_0(\varphi),\zeta\rangle
 =\int_{\Omega_T}q|\partial_tv_\varphi|^{p-2}
                     \partial_tv_\varphi\,w_\zeta\,dx\,dt.
\end{equation}
\end{proposition}
\begin{proof}
Let $u_\delta$ be the unique solution to \eqref{eq:nld-problem} with exterior condition $\delta\varphi$, $0<\delta\leq 1$, and let $v_\varphi$ be the solution of \eqref{eq:v-phi}. Set $z_\delta \vcentcolon =\delta^{-1}u_\delta$, $h_\delta \vcentcolon =z_\delta-v_\varphi$,
$a_\delta \vcentcolon =\partial_tz_\delta$, $b \vcentcolon =\partial_tv_\varphi$, and
$\lambda \vcentcolon =\delta^{p-2}$. The function $h_\delta$ has zero exterior
and initial data and satisfies
\[
 (\partial_t^2+(-\Delta)^r\partial_t +(-\Delta)^s)h_\delta=-\lambda qN_p(a_\delta).
\]
By the forward estimates and Lemma~\ref{lem:nld-embedding}, all the
velocities in this equality lie in $L^p(\Omega_T)$. The difference energy
inequality (see Section~\ref{sec:nonlinear-forward}) is
\begin{equation}\label{eq:nld-difference}
 E_{h_\delta}(t)+\int_0^t \|(-\Delta)^{r/2}\partial_th_\delta\|_{L^2(\R^n)}^2\,d\rho
 +\lambda\int_{\Omega_t}qN_p(a_\delta)(a_\delta-b)\leq0.
\end{equation}
The source product
is integrable by weighted H\"older's inequality. 

Let us introduce the absolutely continuous measure $\mu=q\,dx\,dt$. Nonnegativity of the first two terms and
H\"older's inequality imply
\[
 \int_{\Omega_t}|a_\delta|^p\,d\mu
 \leq\int_{\Omega_t} N_p(a_\delta)b\,d\mu
 \leq\|a_\delta\|_{L^p(\Omega_t,\mu)}^{p-1}\|b\|_{L^p(\Omega_t,\mu)}.
\]
Thus
\begin{equation}\label{eq:nld-weight-bound}
 \begin{split}
     &\|a_\delta\|_{L^p(\Omega_t,\mu)}\leq\|b\|_{L^p(\Omega_t,\mu)},\quad 0<t\leq T,\\
     &\sup_{0\leq t\leq T}E_{h_\delta}(t)
 +\int_0^T \langle (-\Delta)^{r/2}\partial_th_\delta,(-\Delta)^{r/2}\partial_th_\delta\rangle_{L^2(\R^n)}\,dt
 \leq C\lambda\|b\|_{L^p(\Omega_T,\mu)}^p.
 \end{split}
\end{equation}
The second estimate follows from \eqref{eq:nld-difference} by discarding
$\lambda\int_{\Omega_t} q|a_\delta|^p$ and using the first estimate on each time
interval. In particular $a_\delta\to b$ strongly in $L^2(\Omega_T)$.

Since $0\leq q\in L^\infty$ and $|\Omega_T|<\infty$, this convergence
also holds in $\mu$-measure. Boundedness in $L^p(\Omega_T,\mu)$ implies
$a_\delta\rightharpoonup b$ in $L^p(\Omega_T,\mu)$, since every weakly convergent
subsequence has the same limit, which follows by a.e.~convergence. For completeness the uniform monotonicity inequality yields directly
\begin{align*}
 2^{2-p}\|a_\delta-b\|_{L^p(\Omega_T,\mu)}^p
 &\leq\int_{\Omega_T}(N_p(a_\delta)-N_p(b))(a_\delta-b)\,d\mu\\
 &\leq-\int_{\Omega_T} N_p(b)(a_\delta-b)\,d\mu\to 0.
\end{align*}
The second inequality follows from \eqref{eq:nld-difference} with $t=T$.
Here $N_p(b)\in L^{p'}(\Omega_T,\mu)$, so the last limit is weak convergence.
This proves the required weighted strong convergence. Alternatively,
\eqref{eq:nld-weight-bound}, lower semicontinuity, and uniform convexity
give the same conclusion.

The pointwise Nemytskii estimate and weighted H\"older give
\[
 \|N_p(a_\delta)-N_p(b)\|_{L^{p'}(\Omega_T,\mu)}
 \leq C_p(\|a_\delta\|_{L^p(\Omega_T,\mu)}+\|b\|_{L^p(\Omega_T,\mu)})^{p-2}
                  \|a_\delta-b\|_{L^p(\Omega_T,\mu)}\to 0.
\]
The passage to an unweighted source norm is justified by
\begin{equation}\label{eq:nld-weight-source}
 \|q(N_p(a_\delta)-N_p(b))\|_{L^{p'}(\Omega_T)}^{p'}
 \leq\|q\|_\infty^{p'-1}
       \|N_p(a_\delta)-N_p(b)\|_{L^{p'}(\Omega_T,\mu)}^{p'}\to 0.
\end{equation}
If $q=0$ a.e., this assertion and the response formula are immediate.
Weighted convergence alone would not imply unweighted convergence of
$a_\delta$ in $L^p$ on the zero set of $q$; that assertion is not used.

Finally apply Lemma~\ref{lem:nld-green} to
$H=u_\delta-\delta v_\varphi$, whose source is
$f=-qN_p(\partial_tu_\delta)$. Using the definition of the DN map
\eqref{eq:nld-DN} and positive homogeneity for $\delta>0$ gives the exact
identity
\begin{equation}\label{eq:nld-exact-response}
 \delta^{-(p-1)}
 \langle\Lambda_q(\delta\varphi)-\delta\Lambda_0(\varphi),\zeta\rangle
 =\int_{\Omega_T}qN_p(a_\delta)w_\zeta\,dx\,dt.
\end{equation}
By \eqref{eq:nld-weight-source} and $w_\zeta\in L^p(\Omega_T)$,
the right side converges to that in \eqref{eq:nld-response}.
\end{proof}

\subsection{Uniqueness of nonlinear potentials}

Our uniqueness theorem for homogeneous potentials reads as follows.

\begin{theorem}\label{thm:nld-uniqueness}
 Let $\Omega\subset\R^n$ be a bounded Lipschitz domain and $T>0$. Suppose the exponents $0< r<s$ satisfy one of the following conditions:
\begin{enumerate}[(i)]
\item\label{R:cond 1 exp Runge uniqueness}
$s\in\N$, $r\in\R_+\setminus\N$;
\item\label{R:cond 2 exp Runge uniqueness}
$s\in\R_+\setminus\N$, $r\in\N_0$;

\item\label{R:cond 3 exp Runge uniqueness}
Assume $n\geq2$.
$r,s\in\R_+\setminus\N$ and
\begin{equation}
\begin{cases}
s-r\notin\mathbb N,
&
\quad \text{if } n\in2\mathbb N,
\\[1mm]
s-r\notin\frac12\mathbb N,
& 
\quad \text{if } n\in 2\mathbb N+1.
\end{cases}
\end{equation}
\end{enumerate}
Assume that $2\leq p<\infty$ satisfies $2<p\leq 2+4r/n$. Let $q_1,q_2\in L^\infty(\Omega_T)$ be nonnegative. Let $W_1,W_2\subset\Omega_e$ be nonempty open sets. If for every
$\varphi\in C_c^\infty((W_1)_T)$ and all sufficiently small
$\delta>0$,
\[
 \Lambda_{q_1}(\delta\varphi)|_{(W_2)_T}
 =\Lambda_{q_2}(\delta\varphi)|_{(W_2)_T},
\]
then $q_1=q_2$ a.e. in $\Omega_T$.
\end{theorem}
\begin{proof}
Let us write $Q=q_1-q_2$. Subtracting the two proved response limits from Proposition~\ref{prop:nld-response} yields
\begin{equation}\label{eq:nld-integral}
 \int_{\Omega_T}QN_p(\partial_tv_\varphi)w_\zeta\,dx\,dt=0
\end{equation}
for every $\varphi\in C_c^\infty((W_1)_T)$ and
$\zeta\in C_c^\infty((W_2)_T)$. Recall that $w_\zeta$ is the unique solution from the adjoint problem with exterior value $\zeta$, while $v_\varphi$ solves \eqref{eq:v-phi}.\\

\noindent\textit{First limit: adjoint Runge.}

Fix $\varphi$ and $\psi\in C_c^\infty(\Omega_T)$. The adjoint Runge approximation theorem (Theorem~\ref{R:thm: Runge approx in X}) gives
$w_{\zeta_k}-\zeta_k\to\psi$ in $X^1_{s,r,T}(\Omega_T)$, hence
$w_{\zeta_k}|_{\Omega_T}\to\psi$ in $L^p(\Omega_T)$ by Lemma~\ref{lem:nld-embedding}.
Since $\partial_tv_\varphi\in L^p(\Omega_T)$, one has
$QN_p(\partial_tv_\varphi)\in L^{p'}(\Omega_T)$.
H\"older's inequality in \eqref{eq:nld-integral} therefore gives
\begin{equation}\label{eq:nld-smooth-test}
 \int_{\Omega_T}QN_p(\partial_tv_\varphi)\psi\,dx\,dt=0.
\end{equation}

\noindent\textit{Second limit: forward velocity approximation.}

The first limit shows that, for every
$\varphi\in C_c^\infty((W_1)_T)$,
\[
    QN_p(\partial_t v_\varphi)=0
    \quad\text{a.e. in }\Omega_T,
\]
since $QN_p(\partial_t v_\varphi)\in L^{p'}(\Omega_T)$
and \eqref{eq:nld-smooth-test} holds for every
$\psi\in C_c^\infty(\Omega_T)$.

Fix $\eta\in C_c^\infty(\Omega_T)$ and define
\[
    \Phi(x,t)=\int_0^t\eta(x,\rho)\,d\rho.
\]
Then $\Phi\in X^1_{s,r,0}(\Omega_T)$.
By Theorem~\ref{R:thm: Runge approx in X}, there exists a sequence
$(\varphi_k)_{k\in\N}\subset C_c^\infty((W_1)_T)$ such that
\[
    v_{\varphi_k}-\varphi_k\to\Phi
    \quad\text{in }X^1_{s,r,0}(\Omega_T).
\]
In particular,
\[
    \partial_t v_{\varphi_k}\to\eta
    \quad\text{in }L^2(\Omega_T).
\]
Passing to a subsequence, we may assume that
\[
    \partial_t v_{\varphi_k}(x,t)\to\eta(x,t)
    \quad\text{for a.e. }(x,t)\in\Omega_T.
\]
For each $k\in\N$, the first limit gives
\[
    QN_p(\partial_t v_{\varphi_k})=0
    \quad\text{a.e. in }\Omega_T.
\]
Taking the union of the corresponding exceptional null sets
and using the continuity of $N_p$, we obtain
\[
    QN_p(\eta)=0
    \quad\text{a.e. in }\Omega_T.
\]
For any compact set $K\Subset\Omega_T$, choose
$\eta\in C_c^\infty(\Omega_T)$ with $\eta=1$ on $K$.
Since $N_p(1)=1$, this yields $Q=0$ a.e. on $K$.
A countable compact exhaustion of $\Omega_T$ therefore gives
$Q=0$ a.e. in $\Omega_T$, completing the proof.
\end{proof}

\appendix
\section{Smoothing of higher-order heat equations}
\label{R:sec: smoothing higher order heat eq}

\begin{lemma}
\label{R:lemma: smoothing}
    Let $k\in \N$, $\tau>0$, and $u_0 \in L^2(\R^n)$. Define $G_k$ by
    \[
        G_k(\tau,x)=\ifourier_{\xi\to x} \!\left(e^{-\tau|\xi|^{2k}}\right).
    \]
    Then the function $u\in L^2(\R^n)$, defined by
\[
    u(x) \vcentcolon= (G_k(\tau)\ast u_0)(x),\quad x\in\R^n,
\]
is real-analytic in $\R^n$.  More precisely, there exists a constant $C=C(\tau,n,k)>0$ such that for every multiindex $\alpha\in\mathbb{N}_0^n$ we have the following Gevrey-type estimate:
\[
    \|\partial_x^\alpha u\|_{L^\infty(\mathbb{R}^n)}
    \leq C^{|\alpha|+1} (|\alpha|!)^{\frac{1}{2k}} \|u_0\|_{L^2(\mathbb{R}^n)}.
\]
\end{lemma}

\begin{proof}
    Using the Fourier transform, we write
\[
    u(x)
    =
    (2\pi)^{-n}
    \int_{\mathbb{R}^n}
    e^{i x\cdot \xi}
    e^{-\tau |\xi|^{2k}}
    \widehat{u_0}(\xi)\,d\xi.
\]
Differentiating under the integral sign yields, for any multiindex $\alpha$,
\[
    \partial_x^\alpha u(x)
    =
    (2\pi)^{-n}
    \int_{\mathbb{R}^n}
    (i\xi)^\alpha
    e^{i x\cdot \xi}
    e^{-\tau |\xi|^{2k}}
    \widehat{u_0}(\xi)\,d\xi.
\]
Taking absolute values and applying the Cauchy--Schwarz inequality and Plancherel's theorem, we obtain
\begin{equation}
\label{R:eq: L2 estimate poly 1}
\begin{split}
    |\partial_x^\alpha u(x)|
    &\leq
    (2\pi)^{-n}
    \|\widehat{u_0}\|_{L^2(\mathbb R^n)}
    \left\|
        |\xi|^{|\alpha|}e^{-\tau|\xi|^{2k}}
    \right\|_{L^2(\mathbb R^n)}
    \\
    &=
    (2\pi)^{-n/2}
    \|u_0\|_{L^2(\mathbb R^n)}
    \left\|
        |\xi|^{|\alpha|}e^{-\tau|\xi|^{2k}}
    \right\|_{L^2(\mathbb R^n)}.
\end{split}
\end{equation}
Here we used the normalization
\[
    \|\widehat{u_0}\|_{L^2(\mathbb R^n)}
    =(2\pi)^{n/2}\|u_0\|_{L^2(\mathbb R^n)}.
\]
It remains to estimate the weighted exponential factor. We assert that
\begin{equation}
\label{R:eq: L2 estimate poly 2}
    \left\|
        |\xi|^{|\alpha|}
        e^{-\tau |\xi|^{2k}}
    \right\|_{L^2(\mathbb{R}^n)}
    \leq
    C^{|\alpha|+1} (|\alpha|!)^{\frac{1}{2k}}
\end{equation}
for some constant $C=C(k,\tau,n)>0$.
Let us set $m\vcentcolon = |\alpha|$ and $\omega_n = |\partial B_1(0)|$. Then, we use polar coordinates and the change of variables
\[
    s = 2\tau r^{2k},
    \qquad
    r = \left(\frac{s}{2\tau}\right)^{\frac{1}{2k}},
\]
with
\[
    dr = \frac{1}{2k}(2\tau)^{-\frac{1}{2k}} s^{\frac{1}{2k}-1}\,ds,
\]
to find
\[
\begin{aligned}
    \int_{\mathbb{R}^n}
    |\xi|^{2m} e^{-2\tau |\xi|^{2k}}\,d\xi
    &=
    \omega_n \int_0^\infty r^{2m+n-1} e^{-2\tau r^{2k}}\,dr\\
    &=\omega_n  \frac{(2\tau)^{-\frac{2m+n}{2k}}}{2k}  \int_0^\infty    s^{\frac{2m+n}{2k}-1} e^{-s}\,ds\\
    &= C(2\tau)^{-\frac{m}{k}}   \Gamma\!\left(\frac{2m+n}{2k}\right)\\
    &\leq C^{m+1}\Gamma\!\left(\frac{2m+n}{2k}\right)
\end{aligned}
\]
for some $C=C(k,\tau,n)>0$. Therefore,
\[
    \left\|        |\xi|^{m} e^{-\tau |\xi|^{2k}}    \right\|_{L^2(\mathbb{R}^n)}    \leq     C^{m+1}\Gamma\!\left(\frac{2m+n}{2k}\right)^{1/2}.
\]
Next we derive an upper bound for the Gamma function. We start from the classical identity
\[
    \Gamma(x)=\sqrt{2\pi}x^{x-1/2}e^{-x}e^{\mu(x)},\quad x>0,
\]
where $\mu$ satisfies $0<\mu(x)<1/(12x)$. Using
\[
     \frac{2m+n}{2k}\geq \frac{n}{2k},
\]
we deduce that there exists a constant $C=C(k,n)>0$ such that
\[
    \Gamma\left(\frac{2m+n}{2k}\right)\leq C\left(\frac{2m+n}{2k}\right)^{\frac{2m+n}{2k}}e^{-\frac{2m+n}{2k}}.
\]
Since
\[
    \frac{2m+n}{2k}\leq Cm,\quad m\geq 1,
\]
we find
\[
\begin{split}
    \left(\frac{2m+n}{2k}\right)^{\frac{2m+n}{2k}}
    &=e^{\frac{2m+n}{2k}\log \left(\frac{2m+n}{2k}\right)}\leq e^{\frac{2m+n}{2k}\log (Cm)}\\
    &= e^{\frac{2m+n}{2k}\log (C)+\frac{2m+n}{2k}\log (m)}\leq e^{Cm+\left(\frac{m}{k}+C\right)\log (m)}\\
     &\leq e^{Cm+\frac{m}{k}\log (m)}\leq C^{m+1} m^{m/k}
\end{split}
\]
for all $m\geq 1$. Hence
\[
    \Gamma\left(\frac{2m+n}{2k}\right)
    \leq C^{m+1} m^{m/k}e^{-m/k}.
\]
By Stirling's formula,
\[
    m!\sim \sqrt{2\pi m}\,m^m e^{-m},
\]
and therefore
\[
    m^{m/k}e^{-m/k}
    \leq C (m!)^{1/k}
\]
for all sufficiently large $m$. Enlarging $C$ to absorb the finitely many
remaining values of $m$, we obtain
\[
    \Gamma\!\left(\frac{2m+n}{2k}\right)
    \leq C^{m+1}(m!)^{1/k}
    \qquad\text{for all }m\in\mathbb N_0.
\]
Taking square roots yields
\[    
    \Gamma\!\left(\frac{2m+n}{2k}\right)^{1/2}    \leq    C^{m+1} (m!)^{\frac{1}{2k}}.
\]
This establishes formula \eqref{R:eq: L2 estimate poly 2}.

Combining the estimates \eqref{R:eq: L2 estimate poly 1} and \eqref{R:eq: L2 estimate poly 2} yields
\[
    |\partial_x^\alpha u(x)|
    \leq
    C^{|\alpha|+1} (|\alpha|!)^{\frac{1}{2k}} \|u_0\|_{L^2(\mathbb{R}^n)}.
\]
Finally, since
\[
    (|\alpha|!)^{1/(2k)}
    \leq
    |\alpha|!,
\]
the Taylor series of $u$ has positive radius of convergence at every point. Hence $u$ is real-analytic in $\R^n$.
\end{proof}

\medskip
\noindent\textbf{Acknowledgments.} P.~Zimmermann is supported by the Swiss National Science Foundation (SNSF) under Grant Nos. 214500 and 239130. The author is particularly grateful to Katya Krupchyk for numerous stimulating and valuable discussions throughout the development of this work. Part of this research was carried out during a visit to the University of California, Irvine. The author thanks the institution for its hospitality.

\medskip 
\noindent\textbf{Declaration of AI-assisted tools.} 
Following the completion of the initial manuscript, the author used OpenAI's ChatGPT to assist with improvements to the organization and clarity of the exposition and the identification of potential proof gaps or errors. All mathematical statements and proofs were independently reviewed and verified by the author, who takes full responsibility for the content of the manuscript.

\bibliography{refs} 
	
\bibliographystyle{alpha}
    
\end{document}